\documentclass[journal]{IEEEtran}
\usepackage[T1]{fontenc}
\usepackage[latin9]{inputenc}
\usepackage{color}
\usepackage{babel}
\usepackage{mathtools}
\usepackage{amsmath}
\usepackage{amsthm}
\usepackage{amssymb}
\usepackage{graphicx}

\makeatletter
\theoremstyle{plain}
\newtheorem{thm}{\protect\theoremname}
\theoremstyle{remark}
\newtheorem{rem}{\protect\remarkname}
\theoremstyle{definition}
\newtheorem{defn}{\protect\definitionname}
\theoremstyle{plain}
\newtheorem{cor}{\protect\corollaryname}
\theoremstyle{plain}
\newtheorem{prop}{\protect\propositionname}
\theoremstyle{plain}
\newtheorem{fact}{\protect\factname}

\usepackage[mathscr]{eucal}
\usepackage{epsfig,epsf,psfrag}
\usepackage{amssymb,amsmath,amsthm,latexsym}
\usepackage{amsmath,graphicx,xcolor,url}
\usepackage[caption=false]{subfig} 
\usepackage{fixltx2e}
\usepackage{array}
\usepackage{verbatim}
\usepackage{bm}
\usepackage{algorithmic}
\usepackage{algorithm}
\usepackage{verbatim}
\usepackage{textcomp}
\usepackage{mathrsfs,overpic}
\usepackage{epstopdf}
\usepackage{amsfonts}
\usepackage[numbers]{natbib}
\mathtoolsset{showonlyrefs} 

\usepackage{tikz}
\usepackage{float}
\usepackage{tabularx}
\usepackage{multirow}
\usetikzlibrary{patterns}

\newcommand{\Unif}{\mathrm{Unif}}

\def\Ent{\operatorname{Ent}}

\def\conv{{\rm conv}}
\def\conc{{\rm conc}}

\def\1{\mathbf{1}}

\def\d{{\text {\rm d}}}

\catcode`~=11 \def\UrlSpecials{\do\~{\kern -.15em\lower .7ex\hbox{~}\kern .04em}} \catcode`~=13 

\allowdisplaybreaks[1]

\newcommand{\calX}{\mathcal{X}}

\newcommand{\bA}{\mathbf{A}}

\newcommand{\bI}{\mathbf{I}}

\newcommand{\bL}{\mathbf{L}}

\newcommand{\rmd}{\mathrm{d}}

\newcommand{\rme}{\mathrm{e}}

\newcommand{\bbR}{\mathbb{R}}

\DeclareMathAlphabet{\mathbsf}{OT1}{cmss}{bx}{n}
\DeclareMathAlphabet{\mathssf}{OT1}{cmss}{m}{sl}

\DeclareSymbolFont{bsfletters}{OT1}{cmss}{bx}{n}  
\DeclareSymbolFont{ssfletters}{OT1}{cmss}{m}{n}
\DeclareMathSymbol{\bsfGamma}{0}{bsfletters}{'000}
\DeclareMathSymbol{\ssfGamma}{0}{ssfletters}{'000}
\DeclareMathSymbol{\bsfDelta}{0}{bsfletters}{'001}
\DeclareMathSymbol{\ssfDelta}{0}{ssfletters}{'001}
\DeclareMathSymbol{\bsfTheta}{0}{bsfletters}{'002}
\DeclareMathSymbol{\ssfTheta}{0}{ssfletters}{'002}
\DeclareMathSymbol{\bsfLambda}{0}{bsfletters}{'003}
\DeclareMathSymbol{\ssfLambda}{0}{ssfletters}{'003}
\DeclareMathSymbol{\bsfXi}{0}{bsfletters}{'004}
\DeclareMathSymbol{\ssfXi}{0}{ssfletters}{'004}
\DeclareMathSymbol{\bsfPi}{0}{bsfletters}{'005}
\DeclareMathSymbol{\ssfPi}{0}{ssfletters}{'005}
\DeclareMathSymbol{\bsfSigma}{0}{bsfletters}{'006}
\DeclareMathSymbol{\ssfSigma}{0}{ssfletters}{'006}
\DeclareMathSymbol{\bsfUpsilon}{0}{bsfletters}{'007}
\DeclareMathSymbol{\ssfUpsilon}{0}{ssfletters}{'007}
\DeclareMathSymbol{\bsfPhi}{0}{bsfletters}{'010}
\DeclareMathSymbol{\ssfPhi}{0}{ssfletters}{'010}
\DeclareMathSymbol{\bsfPsi}{0}{bsfletters}{'011}
\DeclareMathSymbol{\ssfPsi}{0}{ssfletters}{'011}
\DeclareMathSymbol{\bsfOmega}{0}{bsfletters}{'012}
\DeclareMathSymbol{\ssfOmega}{0}{ssfletters}{'012}

\newcommand{\Bern}{\mathrm{Bern}}

\DeclareMathOperator{\diag}{diag}

\DeclareMathOperator{\supp}{supp}

\newtheorem{definition}{Definition}

\newcommand{\qednew}{\nobreak \ifvmode \relax \else
      \ifdim\lastskip<1.5em \hskip-\lastskip
      \hskip1.5em plus0em minus0.5em \fi \nobreak
      \vrule height0.75em width0.5em depth0.25em\fi}

\usepackage{bm,bbm}

\newcommand{\bzero}{\mathbf{0}}
\newcommand{\bone}{\mathbbm{1}}
\usepackage{babel}
\usepackage{accents} 
\usepackage{enumerate}
\usepackage{verbatim}
\usepackage{mathrsfs}
\usepackage{tikz}
\usepackage{bm,bbm}

\makeatother

\providecommand{\corollaryname}{Corollary}
\providecommand{\definitionname}{Definition}
\providecommand{\factname}{Fact}
\providecommand{\propositionname}{Proposition}
\providecommand{\remarkname}{Remark}
\providecommand{\theoremname}{Theorem}

\begin{document}
\title{R\'enyi--Sobolev Inequalities and  Connections to Spectral Graph Theory}
\author{Lei Yu, \IEEEmembership{Member, IEEE,} and Hao Wu
\thanks{L. Yu and H. Wu are  with the School of Statistics and Data Science,
LPMC, KLMDASR, and LEBPS, Nankai University, Tianjin 300071, China
(e-mails: leiyu@nankai.edu.cn and haowu@mail.nankai.edu.cn). This
work was supported by the National Key Research and Development Program
of China under grant 2023YFA1009604, the NSFC under grant 62101286,
and the Fundamental Research Funds for the Central Universities of
China (Nankai University) under grant 054-63243076. The corresponding
author is L. Yu. } }
\maketitle
\begin{abstract}
In this paper, we generalize the log-Sobolev inequalities to R\'enyi--Sobolev
inequalities by replacing the entropy with  the two-parameter entropy,
which is a generalized version of entropy and closely related to R\'enyi
divergences. We derive the sharp nonlinear dimension-free version
of this kind of inequalities. Interestingly, the resultant inequalities
show a transition phenomenon depending on the parameters. We then
connect R\'enyi--Sobolev inequalities to contractive and data-processing
inequalities, concentration inequalities, and spectral graph theory.
Our proofs in this paper are based on the information-theoretic characterization
of the R\'enyi--Sobolev inequalities, as well as the method of types.
\end{abstract}


\section{Introduction }

Let $\mathcal{X}$ be a finite set and $\bL$ be a $|\calX|\times|\calX|$
 symmetric\footnote{For simplicity, we assume $\bL$ to be symmetric. The results obtained
in this paper can be easily to extended to the asymmetric case. } square matrix such that $L_{x,y}\ge0$ for $x\neq y$ and $\bL\mathbf{1}=\bzero$,
where $\mathbf{1}$ denotes the vector of ones and $\bzero$ denotes
the zero vector. Let $T_{t}:=\rme^{t\bL}$ for $t\ge0$, where $\rme^{\bL}$
denotes the {\em matrix exponential} of $\bL$.  The operator
$T_{t}$ is known as an {\em Markov operator}. In addition, $\{T_{t}\}_{t\ge0}$
forms a {\em Markov semigroup}, since it  satisfies the semigroup
property $T_{t+s}=T_{t}T_{s}=T_{s}T_{t}$ for all $s,t\ge0$. For
more details on Markov operators and Markov semigroups, the reader
is referred to \cite{bakry2004functional,bakry2013analysis,rudnicki2002markov}.

Let $\pi$ be a stationary distribution corresponding to $\{T_{t}\}_{t\ge0}$,
which can be thought as a vector of length $|\calX|$. In other words,
$\pi$ is a distribution satisfying the equation $\pi^{\top}=\pi^{\top}T_{t}$
for all $t\ge0$ or, equivalently, $\pi^{\top}\bL=\bzero^{\top}$.
Without loss of generality, we assume $\pi(x)>0$ for all $x\in\calX$,
since otherwise, we reset  $\calX$ to $\left\{ x:\pi(x)>0\right\} $. 

 As usual, denote $\langle f\rangle:=\langle f\rangle_{\pi}:=\mathbb{E}_{\pi}[f]=\sum_{x}\pi(x)f(x)$
and $\langle f,g\rangle:=\langle fg\rangle$ for two real-valued functions
$f$ and $g$ defined on $\mathcal{X}$. Denote $(\bL f)(x)=\sum_{y}L_{x,y}f(y)$.
Define the carr\'e du champ operator (i.e., squared field operator),
which is the symmetric bilinear map from $\calX\times\calX$ into
$\calX$ defined as 
\begin{align}
\Gamma(f,g) & =\frac{1}{2}(\bL(fg)-f\bL(g)-g\bL(f)),\label{eq:cdcoperator}
\end{align}
which can be expressed as 
\begin{align}
\Gamma(f,g)(x) & =\frac{1}{2}\sum_{y}L_{x,y}(f(y)-f(x))(g(y)-g(x)).\label{eq:positivity}
\end{align}
 Then the {\em Dirichlet form} of $\{T_{t}\}_{t\ge0}$ is defined
as 
\begin{align}
\mathcal{E}(f,g) & :=\langle\Gamma(f,g)\rangle_{\pi}=-\langle\bL f,g\rangle_{\pi}\nonumber \\
 & =-\sum_{x,y}L_{x,y}f(y)g(x)\pi(x).\label{eqn:dirichlet_form}
\end{align}
The {\em normalized Dirichlet form} of $\{T_{t}\}_{t\ge0}$ is
defined as 
\begin{equation}
\overline{\mathcal{E}}(f,g):=\frac{\mathcal{E}(f,g)}{\langle f,g\rangle_{\pi}}.
\end{equation}

We now extend the definitions of the Dirichlet form and its normalized
version to the $n$-dimensional Cartesian product space $\mathcal{X}^{n}$.
Let $T_{t}^{\otimes n}$ be the product semigroup on $\mathcal{X}^{n}$
induced by $T_{t}$. The infinitesimal generator $\bL^{\oplus n}$
of $T_{t}^{\otimes n}$ is given by 
\begin{equation}
\bL^{\oplus n}f(x^{n})=\sum_{i=1}^{n}\bL f(x^{\setminus k},\cdot),\label{eq:-27}
\end{equation}
where $x^{\setminus k}:=(x_{1},\ldots,x_{k-1},x_{k+1},\ldots,x_{n})\in\calX^{n-1}$
denotes the subvector of $x^{n}$ with the $k$-th coordinate removed.
That is, in the summand in \eqref{eq:-27}, the infinitesimal generator
$\bL$ acts on the $k$-th coordinate of $f$ with other coordinates
held fixed. Correspondingly, the ``carr\'e du champ\textquotedblright{}
operator for this case is 
\begin{align}
\Gamma^{\oplus n}(f,g) & =\frac{1}{2}(\bL^{\oplus n}(fg)-f\bL^{\oplus n}(g)-g\bL^{\oplus n}(f))\label{eq:cdcoperator-2}\\
 & =\frac{1}{2}\sum_{y^{n}}\bL_{\cdot,y^{n}}^{\oplus n}(f(y^{n})-f(\cdot))(g(y^{n})-g(\cdot)).\label{eq:cdc}
\end{align}
 For two real-valued functions $f$ and $g$ defined on $\mathcal{X}^{n}$,
let 
\begin{equation}
\psi(x^{\setminus k}):=\mathcal{E}\big(f(x^{\setminus k},\cdot),g(x^{\setminus k},\cdot)\big)
\end{equation}
be the action of the Dirichlet form $\mathcal{E}$ on the $k$-th
coordinates of $f$ and $g$ with other coordinates held fixed. Then,
the Dirichlet form of $f$ and~$g$ and its normalized version are
respectively given by 
\begin{align}
\mathcal{E}_{n}(f,g) & :=\langle\Gamma^{\oplus n}(f,g)\rangle_{\pi^{n}}=-\langle\bL^{\oplus n}f,g\rangle_{\pi^{n}}\nonumber \\
 & =\sum_{k=1}^{n}\sum_{x^{\setminus k}\in\mathcal{X}^{n-1}}\psi(x^{\setminus k})\pi^{n-1}(x^{\setminus k})\\
 & =\sum_{k=1}^{n}\mathbb{E}\left[\psi(X^{\setminus k})\right],\label{eq:FImaten-1}\\
\overline{\mathcal{E}}_{n}(f,g) & =\frac{\mathcal{E}_{n}(f,g)}{\langle f,g\rangle_{\pi^{n}}}.\label{eqn:norm_df-1}
\end{align}

In the following, we focus on the normalized Dirichlet form $\overline{\mathcal{E}}_{n}(f,f^{q-1})$
with $q>1$. The motivation of considering this kind of Dirichlet
forms comes from the following argument. Define the function $\phi:[0,\infty)\to\bbR$
as 
\begin{equation}
\phi(t):=-\frac{1}{n}\ln\|T_{t}f\|_{q}\,,
\end{equation}
where $\|g\|_{q}:=\left(\mathbb{E}_{\pi}[g^{q}]\right)^{1/q}$ for
nonnegative function $g$. Then, one can check by direct differentiation
that 
\begin{align}
\phi'(0) & =\left.\frac{1}{n}\overline{\mathcal{E}}_{n}\big(T_{t}f,(T_{t}f)^{q-1}\big)\right|_{t=0}=\frac{1}{n}\overline{\mathcal{E}}_{n}\big(f,f^{q-1}\big).\label{eq:phider}
\end{align}

In addition to the Dirichlet form, the other quantity involved in
log-Sobolev inequalities is the entropy of a nonnegative function~$f$. 

\begin{definition} For a nonnegative function $f$, the entropy and
the normalized entropy of $f$ are respectively defined as 
\begin{align}
\Ent(f) & :=\mathbb{E}_{\pi}[f\ln f]-\mathbb{E}_{\pi}[f]\ln\mathbb{E}_{\pi}[f],\nonumber \\
\overline{\Ent}(f) & :=\frac{\Ent(f)}{\mathbb{E}_{\pi}[f]}.\label{eq:FInormalizedentropy-1}
\end{align}
\end{definition} 

Note that these notions of entropy and normalized entropy as defined
here are commonly encountered in functional analysis; see, for example,
\cite{Ledoux_book}. They are related to, but not the same as the
Shannon entropy in classical information theory. Indeed, they bear
more similarity to the relative entropy, in the sense that if $f$
is the Radon--Nikodym derivative ${\rmd Q}/{\rmd\pi}$ of a distribution
$Q$ with respect to $\pi$ (i.e., the function $x\in\calX\mapsto{Q(x)}/{\pi(x)}$
for the finite alphabet case, which is also denoted as $Q/\pi$),
then the entropy (and also the normalized entropy) of $f$ turns into
the relative entropy of $Q$ from $\pi$, i.e., 
\[
D(Q\|\pi):=\sum_{x}Q(x)\ln\frac{Q(x)}{\pi(x)}.
\]
By Jensen's inequality, both the entropy and the normalized entropy
are nonnegative.

\subsection{Log-Sobolev Inequalities and Log-Sobolev Functions}

The log-Sobolev inequalities connect the Dirichlet forms and the entropies
in the following way. 

\begin{definition} \cite{gross1975logarithmic,mossel2013reverse}
For $q\in\mathbb{R}\backslash\{0,1\}$, the {\em $q$-log-Sobolev
inequality with constant $C$} is 
\begin{equation}
\Ent(f^{q})\le C\,\frac{q^{2}}{q-1}\,\mathcal{E}(f,f^{q-1})\label{eq:FIpLSI}
\end{equation}
for positive\footnote{By the monotone convergence theorem, the inequality in \eqref{eq:FIpLSI}
still holds for nonnegative $f$ if $q>1$.}~$f$. 
For $q=0$, the {\em $0$-log-Sobolev inequality with constant $C$}
for positive~$f$ is 
\begin{equation}
\frac{1}{2}\mathrm{Var}(\ln f)\le-C\,\mathcal{E}(f,\frac{1}{f}).\label{eq:FI0LSI}
\end{equation}
For $q=1$, the {\em $1$-log-Sobolev inequality with constant
$C$} for positive~$f$ is 
\begin{equation}
\Ent(f)\le C\,\mathcal{E}(f,\ln f).\label{eq:FI1LSI}
\end{equation}
\end{definition}

The inequalities in \eqref{eq:FI1LSI} and \eqref{eq:FI0LSI} are
limiting versions of the one in \eqref{eq:FIpLSI}, since for positive~$f$,
$\frac{1}{q^{2}}\Ent(f^{q})\to\frac{1}{2}\mathrm{Var}(\ln f)$ as
$q\to0$ and $\frac{1}{q-1}\,\mathcal{E}(f,f^{q-1})\to\mathcal{E}(f,\ln f)$
as $q\to1$. 

Inspired by the log-Sobolev inequality, one can define the log-Sobolev
function. 

\begin{definition} The {\em  log-Sobolev function} for $q\in\mathbb{R}\backslash\{0,1\}$
is defined as
\begin{equation}
\Xi_{q}(\alpha):=\inf_{f:\overline{\Ent}(f^{q})\ge\alpha}\frac{1}{q-1}\overline{\mathcal{E}}(f,f^{q-1}).\label{eq:FIpLSI-1}
\end{equation}
For $q\in\{0,1\}$, 
\begin{align}
\Xi_{1}(\alpha) & :=\inf_{f:\overline{\Ent}(f)\ge\alpha}\frac{\mathcal{E}(f,\ln f)}{\mathbb{E}_{\pi}[f]},\label{eq:FIpLSI-1-3}\\
\Xi_{0}(\alpha) & :=\inf_{f:\frac{1}{2}\mathrm{Var}(\ln f)\ge\alpha}-\overline{\mathcal{E}}(f,f^{-1}).\label{eq:-15}
\end{align}
\end{definition} 

This function characterizes the nonlinear tradeoff between the entropy
and the Dirichlet form. Here, $\Xi_{1}$ is defined by the continuous
extension. The direct extension of the definition of $\Xi_{q}(\alpha)$
to the case $q=0$ does not make sense. So, inspired by \eqref{eq:FI0LSI},
instead we define $\Xi_{0}(\alpha)$ in \eqref{eq:-15} as $\Xi_{q}(q^{2}\alpha)$
with $q=0$ where $\frac{1}{q^{2}}\overline{\Ent}(f^{q})$ for $q=0$
is understood as $\frac{1}{2}\mathrm{Var}(\ln f)$.\footnote{A more consistent way to define the log-Sobolev function is to replace
the constraint in \eqref{eq:FIpLSI-1} with $\overline{\Ent}(f^{q})\ge q^{2}\alpha$.
In other words, the new log-Sobolev function is $\alpha\mapsto\Xi_{q}(q^{2}\alpha)$
with $\Xi_{q}$ denoting the original one defined in \eqref{eq:FIpLSI-1}.
However, in this definition, the effective domains of this function
are not consistent anymore for different $q$'s, since the possible
range of $\overline{\Ent}(f^{q})$ is $[0,\ln\frac{1}{\min_{x}\pi(x)}]$
for any $q>0$ and $[0,\ln\frac{1}{\min_{x}\pi(x)})$ for any $q<0$.} 

By definition, for $q\neq0$, the $q$-log-Sobolev inequality with
constant $C$ holds if and only if $\Xi_{q}(\alpha)\ge\frac{\alpha}{q^{2}C}$.
In other words, for $q\neq0$, the optimal $q$-log-Sobolev constant
$C_{q}^{*}$ is $\sup_{\alpha>0}\frac{\alpha}{q^{2}\Xi_{q}(\alpha)}$.
Similarly, the optimal $0$-log-Sobolev constant $C_{0}^{*}$ is $\sup_{\alpha>0}\frac{\alpha}{\Xi_{0}(\alpha)}$. 

Extending the definition of $\Xi_{q}(\alpha)$ from $T_{t}$ to $T_{t}^{\otimes n}$,
we define for $q\in\mathbb{R}\backslash\{0,1\}$,
\begin{equation}
\Xi_{q}^{(n)}(\alpha):=\inf_{f:\frac{1}{n}\overline{\Ent}(f^{q})\ge\alpha}\frac{1}{(q-1)n}\overline{\mathcal{E}}_{n}(f,f^{q-1}),\label{eqn:def_Xin-1}
\end{equation}
where $f$ is defined on $\mathcal{X}^{n}$. For $q\in\{0,1\}$, extend
the definition of $\Xi_{q}(\alpha)$ from $T_{t}$ to $T_{t}^{\otimes n}$
in a similar way. Since the $q$-log-Sobolev inequality satisfies
the tensorization property (i.e., the $q$-log-Sobolev inequality
is satisfied for $T_{t}$ with constant $C$  if and only if satisfied
for $T_{t}^{\otimes n}$ with the same constant), $\Xi_{q}^{(n)}(\alpha)\ge\frac{\alpha}{C}$
still holds when the $q$-log-Sobolev inequality with constant $C$
holds for $T_{t}$. However, the log-Sobolev function $\Xi_{q}^{(n)}$
itself does not satisfy the tensorization property, but it indeed
satisfies the following quasi-tensorization property. 
\begin{thm}[Quasi-Tensorization \cite{polyanskiy2019improved}]
 \label{thm:strongSLI} For $q\in\mathbb{R}$,  it holds that 
\begin{equation}
\conv\,\Xi_{q}(\alpha)\le\Xi_{q}^{(n)}(\alpha)\le\Xi_{q}(\alpha),\label{eq:FI-strongerlogSobolev}
\end{equation}
where $\conv\,\Xi_{q}$ denotes the lower convex envelope of the function
$\Xi_{q}$. Moreover, this lower bound is asymptotically tight as
$n\to\infty$, which means that 
\begin{equation}
\lim_{n\to\infty}\Xi_{q}^{(n)}(\alpha)=\conv\,\Xi_{q}(\alpha).\label{eqn:asympt_tight_U}
\end{equation}
If additionally, $\Xi_{q}$ is convex, then the lower bound is tight
for any finite $n\ge1$.  In fact, there are two distributions $(Q,R)$
and a number $\lambda\in[0,1]$ such that the sequence of functions
\[
f_{n}(x^{n})=\frac{Q^{k}(x^{k})}{\pi^{k}(x^{k})}\times\frac{R^{n-k}(x_{k+1}^{n})}{\pi^{n-k}(x_{k+1}^{n})}
\]
with $k=\left\lfloor \lambda n\right\rfloor $ satisfies $\frac{1}{n}\overline{\Ent}(f^{q})\ge\alpha$
and $\lim_{n\to\infty}\frac{1}{(q-1)n}\overline{\mathcal{E}}_{n}(f_{n},f_{n}^{q-1})=\conv\,\Xi_{q}(\alpha)$.
Here, $\left\lfloor x\right\rfloor $ denotes the maximum integer
no larger than $x$. 
\end{thm}
\begin{rem}
\label{rem:The-inequality-in-1}The inequality in \eqref{eq:FI-strongerlogSobolev}
can be easily extended to semigroups defined on arbitrary measurable
spaces, since the proof given in \cite{polyanskiy2019improved} only
relies on the chain rule of relative entropies and the fact that conditioning
increases relative entropies. 
\end{rem}
\begin{rem}
\label{rem:The-log-Sobolev-function}The log-Sobolev function $\Xi_{q}$
is in general not convex. This was observed by Gu and Polyanskiy in
\cite[Proposition 26]{gu2023non} for the Potts semigroup. 
\end{rem}
This theorem is a strengthening of the (linear) $q$-log-Sobolev inequality
in \eqref{eq:FIpLSI}. Since in general, the function $\conv\,\Xi_{q}$
is nonlinear, the inequality in \eqref{eq:FI-strongerlogSobolev}
is called the nonlinear $q$-log-Sobolev inequality. To see the relation
between the linear and nonlinear $q$-log-Sobolev inequalities, one
can easily find that the optimal constant $C=C_{q}^{*}$ in the (linear)
$q$-log-Sobolev inequality in \eqref{eq:FIpLSI} satisfies that $\frac{1}{q^{2}C_{q}^{*}}$
for $q\neq0$ and $\frac{1}{C_{0}^{*}}$ for $q=0$ are exactly 
the right-derivatives of $\conv\,\Xi_{q}(\alpha)$ at $\alpha=0$.

\subsection{R\'enyi--Sobolev Functions }

In this paper, we aim at generalizing Polyanskiy and Samorodnitsky's
result by replacing the entropy with a general form, known as the
two-parameter entropy. Given orders $p,q>0$ such that $p\neq q$
and a nonnegative $f$, define the $(p,q)$-entropy of $f$ \cite{yu2021strong_article}
as 
\begin{equation}
\Ent_{p,q}(f):=\frac{pq}{p-q}\ln\frac{\Vert f\Vert_{p}}{\Vert f\Vert_{q}}.\label{eq:ent}
\end{equation}
The $(p,q)$-entropies for other pairs of $(p,q)$ are defined by
continuous extensions. In particular, for $q>0$, define 
\[
\Ent_{q}(f):=\Ent_{q,q}(f):=\lim_{p\uparrow q}\Ent_{p,q}(f)=\overline{\Ent}(f^{q}).
\]
For $q\ge0$ and any nonnegative $f$, the $(0,q)$-entropy and
the $(q,0)$-entropy of $f$ are defined by 
\begin{equation}
\Ent_{0,q}(f),\;\Ent_{q,0}(f):=\Ent_{0}(f):=-\ln\pi(f>0).\label{eq:ent0}
\end{equation}
For $q\ge0$ and any nonnegative $f$, the $(\infty,q)$-entropy and
the $(q,\infty)$-entropy of $f$ are defined by 
\begin{align}
\Ent_{\infty,q}(f),\;\Ent_{q,\infty}(f) & :=\Ent_{\infty}(f):=q\ln\frac{\Vert f\Vert_{\infty}}{\Vert f\Vert_{q}}\nonumber \\
 & =-\ln\mathrm{esssup}_{\pi}(\frac{f^{q}}{\mathbb{E}[f^{q}]}),\label{eq:ent0-1}
\end{align}
where $\mathrm{esssup}_{\pi}(g)$ denotes the essential supremum of
$g$. The two-parameter entropy is nondecreasing in its one parameter
given the other one \cite{yu2021strong_article}. The two-parameter
entropy is a natural generalization of the entropy. 

The two-parameter entropy is closely related to the well-known R\'enyi
divergence. The R\'enyi divergence of order $\gamma\in\mathbb{R}\backslash\{0,1\}$
is defined as 
\[
D_{\gamma}(Q\|\pi):=\frac{1}{\gamma-1}\ln\sum_{x}Q(x)\Big(\frac{Q(x)}{\pi(x)}\Big)^{\gamma-1}.
\]
For $\gamma\in\{0,1,\pm\infty\}$, the R\'enyi divergence of order $\gamma$
is defined by continuous extension. In particular, for $\gamma=1$,
the R\'enyi divergence reduces to the relative entropy. For $\gamma=0$,
$D_{0}(Q\|\pi):=-\ln\pi(Q>0)$, and for $\gamma=0$, $D_{\infty}(Q\|\pi):=-\ln\mathrm{esssup}_{Q}(\frac{Q}{\pi})$. 

By definition, the two-parameter entropy and the R\'enyi divergence
 admit the following intimate relationship: For $q\in\mathbb{R}\backslash\{0\}$
and any $f$ such that $0<\Vert f\Vert_{q}<\infty$, it holds that
\begin{equation}
\Ent_{p,q}(f)=D_{p/q}(Q\|\pi).\label{eq:-3}
\end{equation}
where $Q$ is a distribution such that $\frac{f^{q}}{\mathbb{E}[f^{q}]}=\frac{Q}{\pi}$
with $\frac{Q}{\pi}$ denoting the Radon--Nikodym derivative.

The two-parameter entropy was previously used in \cite{samorodnitsky2008modified_article,polyanskiy2019improved,kirshner2021moment,yu2024graphs,yu2021strong_article}
to investigate the nonlinear versions of log-Sobolev, hypercontractivity,
and Brascamp--Lieb inequalities. 

\begin{definition} Define the {\em  R\'enyi--Sobolev function}
as 
\[
\Xi_{p,q}(\alpha):=\inf_{f:\Ent_{p,q}(f)\ge\alpha}\frac{1}{q-1}\overline{\mathcal{E}}(f,f^{q-1}).
\]
\end{definition} 

Extending the definition of $\Xi_{p,q}(\alpha)$ from $T_{t}$ to
$T_{t}^{\otimes n}$, we define the $n$-dimensional version as 
\begin{equation}
\Xi_{p,q}^{(n)}(\alpha):=\inf_{f:\frac{1}{n}\Ent_{p,q}(f)\ge\alpha}\frac{1}{(q-1)n}\overline{\mathcal{E}}_{n}(f,f^{q-1}).\label{eq:Xi}
\end{equation}
When specialized to the case $p=q$, $\Xi_{p,q}$ and $\Xi_{p,q}^{(n)}$
respectively reduce to $\Xi_{q}$ and $\Xi_{q}^{(n)}$ defined in
\eqref{eq:FIpLSI-1} and \eqref{eqn:def_Xin-1}. By the monotonicity
of the two-parameter entropy in its argument, $\Xi_{p,q}^{(n)}(\alpha)$
is nonincreasing in $p$ given $n,q,\alpha$. 

\subsection{Our Contributions}

Our contributions in this paper are as follows. 
\begin{enumerate}
\item In this paper, we investigate the asymptotics of the R\'enyi--Sobolev
function, and derive the sharp dimension-free bound on the R\'enyi--Sobolev
function. The resultant inequalities are called R\'enyi--Sobolev inequalities.
Interestingly, R\'enyi--Sobolev inequalities show a transition phenomenon
depending on the parameters $(p,q)$. In particular, the R\'enyi--Sobolev
function converges to $\conv\,\Xi_{q}(\alpha)$ as $n\to\infty$ when
$p\le q$, and converges to $0$ when $p>q$. Our proofs are based
on the information-theoretic characterization of the R\'enyi--Sobolev
inequalities, as well as the method of types. 
\item We then connect R\'enyi--Sobolev inequalities to several other topics,
including contractive and data-processing inequalities, concentration
inequalities, and spectral graph theory. In the aspect of contractive
and data-processing inequalities, our R\'enyi--Sobolev inequalities
characterize the instantaneous change of the $q$-norm of $T_{t}f$
for a ``non-flat'' function $f$.  In the aspect of concentration
inequalities, our R\'enyi--Sobolev inequalities are coupled with the
Herbst argument to yield new concentration inequalities, which are
stronger than the ones derived by the Herbst argument coupled with
standard log-Sobolev inequalities. In the aspect of spectral graph
theory, our R\'enyi--Sobolev inequalities are used to characterize
the asymptotics of the spectral radius of the $n$-fold Cartesian
product of a regular graph with itself.  This generalizes the existing
result in \cite{friedman2005generalized,bollobas2018eigenvalues}
from the complete graph $K_{2}$ to arbitrary regular graph. 
\end{enumerate}

\subsection{Organization}

This paper is organized as follows. In Section \ref{sec:Main-Result},
we introduce our main results on R\'enyi--Sobolev inequalities. In
Sections \ref{sec:Connections-to-Contractive}-\ref{sec:Connections-to-Spectral},
we connect R\'enyi--Sobolev inequalities respectively to contractive
and data-processing inequalities, concentration inequalities, and
spectral graph theory. In the remaining part, we provide proofs for
our results. 

\section{\protect\label{sec:Main-Result}Main Result}

\subsection{R\'enyi--Sobolev Inequalities}

Our main result is the following theorem which characterizes the
asymptotics of the R\'enyi--Sobolev function. The proof is deferred
to Section \ref{sec:Proof-of-Theorem}.
\begin{thm}[R\'enyi--Sobolev Inequalities]
\label{thm:RSI}  For $p,q\ge0$,  it holds that 
\begin{equation}
\Xi_{p,q}^{(n)}(\alpha)\ge\phi_{p,q}(\alpha),\label{eq:RenyiSobolev}
\end{equation}
where 
\begin{equation}
\phi_{p,q}(\alpha):=\begin{cases}
\conv\,\Xi_{q}(\alpha), & p\le q\\
0, & p>q
\end{cases}.\label{eq:FI-strongerlogSobolev-1-1}
\end{equation}
Moreover, for $p\ge0,q\ge1$, and $\alpha\in[0,-\ln\min_{x}\pi(x))$,
this lower bound is asymptotically tight as $n\to\infty$, i.e.,
\begin{equation}
\lim_{n\to\infty}\Xi_{p,q}^{(n)}(\alpha)=\phi_{p,q}(\alpha).\label{eq:asympt_tight}
\end{equation}
In fact, for $p\le q$, there are two distributions $(Q,R)$ and
a number $\lambda\in[0,1]$ such that the sequence of functions 
\[
f_{n}(x^{n})=\frac{Q^{k}(x^{k}|\mathcal{T}_{\epsilon}^{(k)}(Q))}{\pi^{k}(x^{k})}\times\frac{R^{n-k}(x_{k+1}^{n}|\mathcal{T}_{\epsilon}^{(n-k)}(R))}{\pi^{n-k}(x_{k+1}^{n})}
\]
with $k=\left\lfloor \lambda n\right\rfloor $ satisfies that $\frac{1}{n}\Ent_{p,q}(f_{n})\ge\alpha$
and $\lim_{n\to\infty}\frac{1}{(q-1)n}\overline{\mathcal{E}}_{n}(f_{n},f_{n}^{q-1})$
can be arbitrarily close to the lower bound $\conv\,\Xi_{q}(\alpha)$
as $\epsilon\downarrow0$.  Here, the $\epsilon$-typical set of
a distribution $Q$  is defined as 
\[
\mathcal{T}_{\epsilon}^{(n)}(Q):=\{x^{n}:|T_{x^{n}}(a)-Q(a)|\le\epsilon Q(a),\forall a\in\calX\},
\]
where $T_{x^{n}}$ is the empirical measure of $x^{n}$. 
\end{thm}
\begin{rem}
\label{rem:The-inequality-in} Similarly to Remark \ref{rem:The-inequality-in-1},
the inequality in \eqref{eq:RenyiSobolev} can be easily extended
to semigroups defined on arbitrary measurable spaces, since the essential
part of the proof given in \cite{polyanskiy2019improved} only relies
on the chain rule of relative entropies and the fact that conditioning
increases relative entropies. Our main contribution in the theorem
above is the equation \eqref{eq:asympt_tight}, i.e., the asymptotic
tightness of the inequality in \eqref{eq:RenyiSobolev}. The proof
for this part relies on the method of types, and hence, cannot be
applied directly to semigroups defined on arbitrary measurable spaces.
\end{rem}
\begin{rem}
The extremal functions in Theorem \ref{thm:RSI} for $p\le q$ are
formed by multiplying two conditional Radon--Nikodym derivatives
on typical sets. Note that for $p=q$, Theorem \ref{thm:strongSLI}
provides another kind of extremal functions which are formed by multiplying
two unconditional Radon--Nikodym derivatives. 
\end{rem}
We now provide more remarks on this theorem.  The range $\alpha\in[0,-\ln\min_{x}\pi(x))$
considered for the equation \eqref{eq:asympt_tight} is reasonable,
since the possible range of $\Ent_{p,q}(f)$ is $[0,-\ln\min_{x}\pi(x)]$.
We exclude the right end point since it corresponds to a trivial case.
In Theorem \ref{thm:RSI}, we prove the asymptotic tightness of the
inequality in \eqref{eq:RenyiSobolev} only for the case $q\ge1$,
which is the most interesting case. It is also interesting to investigate
the case $q<1$. This is left as a future work. 

Some related works include \cite{chafai2004entropies,boucheron2013concentration,raginsky2016strong}.
In these works, the log-Sobolev inequalities are extended to the $\Phi$-Sobolev
inequalities by replacing the entropy with $\Phi$-entropies. Our
extension here is in a similar way, but with the entropy replaced
by R\'enyi entropies. It is known that R\'enyi entropies are not covered
by $\Phi$-entropies. So, R\'enyi--Sobolev inequalities introduced
here are new. Furthermore, only the linear tradeoff between the Dirichlet
form and the $\Phi$-entropy was investigated in \cite{chafai2004entropies,boucheron2013concentration,raginsky2016strong},
but here we investigate a more elaborate tradeoff---the nonlinear
tradeoff between the Dirichlet form and the R\'enyi entropy. When specialized
to the case $p=q$, Theorem \ref{thm:RSI} reduces to Polyanskiy and
Samorodnitsky's nonlinear log-Sobolev inequalities \cite{polyanskiy2019improved}
given in Theorem \ref{thm:strongSLI} for the case $q>1$.

\subsection{\protect\label{subsec:Binary-Example}Binary Example }

Before providing various interpretations of Theorem \ref{thm:RSI},
we first focus on the binary case, and provide explicit expressions
of the nonlinear $q$-log-Sobolev inequalities particularized for
this case. 


Consider $\mathcal{X}=\{0,1\}$, $\pi=\Bern(1/2)$ and $T_{t}$ is
the Bonami--Beckner operator  given by 
\begin{equation}
T_{t}f(y)=f(y)\frac{1+\rme^{-t}}{2}+f(1-y)\frac{1-\rme^{-t}}{2},\;y\in\{0,1\}.\label{eq:FIhypercube}
\end{equation}
 From \eqref{eq:FIhypercube}, we know that the infinitesimal generator
$\bL$ of $T_{t}$ is given by $L_{x,y}=\bone\{x\neq y\}-1/2$ which
is obtained by differentiating $T_{t}$ with respect to $t$ and evaluating
the derivative at $t=0$. Moreover, the Dirichlet form for this case
is 
\begin{align}
\hspace{-0.25in}{\mathcal{E}}_{n}(f,g) & =-\frac{1}{2}\langle\Delta f,g\rangle\quad\mbox{and}\label{eq:FIdiri_cube}\\
\hspace{-0.25in}{\mathcal{E}}_{n}(f,f) & =\frac{2^{-n}}{4}\sum_{(x^{n},y^{n}):x^{n}\sim y^{n}}\big(f(x^{n})-f(y^{n})\big)^{2}\,,\label{eq:Eff}
\end{align}
where $\Delta f(x^{n}):=\sum_{y^{n}:y^{n}\sim x^{n}}(f(y^{n})-f(x^{n}))$,
and $x^{n}\sim y^{n}$ means that $x^{n},y^{n}\in\{0,1\}^{n}$ differ
in exactly one coordinate. 

\begin{figure}

\centering{}\includegraphics[width=0.9\columnwidth]{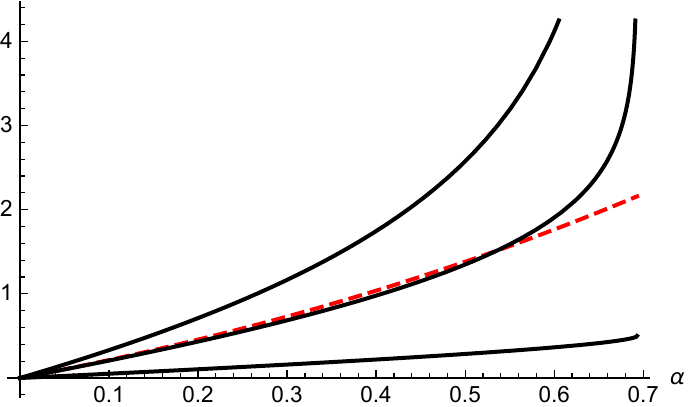}\caption{\protect\label{fig:Illustration-of-the}Illustration of the $q$-log-Sobolev
function $\Xi_{q}$ for $q=0.8,1$, and $2$ (solid curves from top
to bottom) and $\Xi_{q}$ for $q=0$ (dashed curve). }
\end{figure}

For the binary case,  the optimal constant $C$ in the (linear) $q$-log-Sobolev
inequality is $2$; see \cite{gross1975logarithmic}.  In contrast,
the optimal nonlinear $q$-log-Sobolev inequality in Theorem \ref{thm:strongSLI}
was derived by Polyanskiy and Samorodnitsky \cite{polyanskiy2019improved},
showing that 
\begin{align}
\Xi_{q}(\alpha) & =\begin{cases}
{\displaystyle \frac{1}{2(q-1)}\Big(1\!-\!y^{\frac{1}{q}}(1-y)^{\frac{1}{q'}}\!}\\
\qquad{\displaystyle -\!y^{\frac{1}{q'}}(1-y)^{\frac{1}{q}}\Big)}, & q\ne0,1\\
{\displaystyle \Big(\frac{1}{2}-y\Big)\ln\frac{1-y}{y}}, & q=1\\
\frac{1}{4}\left(\rme^{2\sqrt{2\alpha}}+\rme^{-2\sqrt{2\alpha}}\right)-\frac{1}{2}, & q=0
\end{cases},\label{eq:-9}
\end{align}
where $q':=\frac{q}{q-1}$ is the H\"older conjugate of $q$, $y(\alpha):=h^{-1}(\ln2-\alpha)$,
and $h^{-1}:[0,\ln2]\to[0,1/2]$ is the inverse of the binary entropy
function~$h(s)=-s\ln s-(1-s)\ln(1-s)$ with base $\rme$ when its
domain is restricted to $[0,1/2]$. Moreover, they showed that $\Xi_{q}$
is nonnegative, increasing, and convex. The function $\Xi_{q}$ is
plotted in Fig. \ref{fig:Illustration-of-the}. 

Substituting \eqref{eq:-9} to Theorem \ref{thm:RSI} yields the
explicit expressions for R\'enyi--Sobolev inequalities. That is, for
$p,q\ge0$, 
\begin{equation}
\Xi_{p,q}^{(n)}(\alpha)\ge\begin{cases}
\Xi_{q}(\alpha), & p\le q\\
0, & p>q
\end{cases},\label{eq:RenyiSobolev-1}
\end{equation}
which is asymptotically tight for $p\ge0,q\ge1$, and $\alpha\in[0,\ln2)$.

\section{\protect\label{sec:Connections-to-Contractive}Connections to Contractive
and Data-Processing Inequalities}

We now provide two straightforward interpretations for Theorem \ref{thm:RSI}
in the language of contractive  and data-processing inequalities. 

\textbf{Interpretation 1: }Note that $\Ent_{p,q}(f)=0$ if and only
if $f$ is constant. In other words, the two-parameter entropy quantifies
how dramatically the value of the function varies, and hence, can
be seen as a ``distance'' between the function $f$ and the constant
functions. For a function $f$ whose values vary dramatically, specifically,
satisfying $\Ent_{p,q}(f)\ge n\alpha$, how does the $q$-norm of
$T_{t}f$ changes instantaneously at $t=0$?  By \eqref{eq:phider},
Theorem \ref{thm:RSI} provides the sharp dimension-free bound for
this question, which shows a transition phenomenon. Specifically,
if $p\le q$, then the constraint $\Ent_{p,q}(f)\ge n\alpha$ does
not differ too much from $\Ent_{q}(f)\ge n\alpha$ in the asymptotic
regime as $n\to\infty$; on the other hand, if $p>q$, then the constraint
$\Ent_{p,q}(f)\ge n\alpha$ has no effect on the change of the $q$-norm
of $T_{t}f$ in the asymptotic regime as $n\to\infty$. 

\textbf{Interpretation 2: }Note that $T_{t}$ corresponds to a transition
probability matrix, denoted by $P_{X|Y}^{(t)}$. The transition probability
matrix $P_{Y|X}^{(t)}$ is specified by $P_{Y|X}^{(t)}(y|x)=\frac{P_{X|Y}^{(t)}(x|y)\pi(y)}{\pi(x)}$.
That is, $P_{Y|X}^{(t)}=\rme^{t[\pi]^{-1}\bL^{\top}[\pi]}$ where
$[\pi]:=\diag\{\pi(x)\}_{x\in\calX}$. Let $p=1$ and $f=\frac{Q_{X}}{\pi}.$
Then, $\Ent_{p,q}(f)=D_{q}(Q_{X}\|\pi)$. Denote $Q_{Y}^{(t)}$ as
the marginal of $Q_{X}P_{Y|X}^{(t)}.$ Then, $q'\ln\|T_{t}f\|_{q}=D_{q}(Q_{Y}^{(t)}\|\pi)$.
Given the channel $P_{Y|X}^{(t)}$ with stationary distribution $\pi$,
if the input $Q_{X}$ satisfies $D_{q}(Q_{X}\|\pi)\ge\alpha$, how
does $D_{q}(Q_{Y}^{(t)}\|\pi)$ vary instantaneously at $t=0$? Theorem
\ref{thm:RSI} provides the sharp dimension-free bound for this question.

In the interpretations above, the instantaneous changes of $\|T_{t}f\|_{q}$
and $D_{q}(Q_{Y}^{(t)}\|\pi)$ at $t=0$ are considered. In fact,
the quantities $\|T_{t}f\|_{q}$ and $D_{q}(Q_{Y}^{(t)}\|\pi)$ at
a fixed epoch $t>0$ were considered in \cite{yu2021strong_article}
for which the sharp dimension-free bounds were provided. Although
Theorem \ref{thm:RSI} in this paper and the results in \cite{yu2021strong_article}
are closely related, one does not imply the other. On one hand, integrating
the result in Theorem \ref{thm:RSI} with respect to $t$ will yield
a bound on $\|T_{t}f\|_{q}$ and $D_{q}(Q_{Y}^{(t)}\|\pi)$ for a
fixed epoch $t>0$, just as done in \cite{polyanskiy2019improved}.
However, the resulting bound is not asymptotically sharp, and hence,
strictly weaker than the ones in \cite{yu2021strong_article}. On
the other hand, intuitively, differentiating the bounds in \cite{yu2021strong_article}
with respect to $t$ will yield bounds on the change rates of $\|T_{t}f\|_{q}$
and $D_{q}(Q_{Y}^{(t)}\|\pi)$ at $t=0$. However, to this end, we
need to swap optimization operators with the differential operator,
which is not feasible in general.

\section{\protect\label{subsec:Connections-to-Concentration}Connections to
Concentration Inequalities}

We now connect R\'enyi-Sobolev inequalities to concentration of measure.
Instead of the original R\'enyi-Sobolev function defined in \eqref{eq:Xi},
here we consider a generalized version of the R\'enyi-Sobolev function,
and then use this generalized R\'enyi-Sobolev function to strengthen
existing concentration inequalities. In this section, we assume that
the semigroups mentioned here could be defined on arbitrary measurable
spaces, not restricted to finite spaces. 

\subsection{Generalized R\'enyi-Sobolev Function}

Define the ``inverse'' of the R\'enyi-Sobolev function $\Xi_{p,q}^{(n)}$
in \eqref{eq:Xi} as 
\begin{align*}
\Xi_{p,q}^{(n),-1}(t) & :=\sup_{f:\frac{1}{(q-1)n}\overline{\mathcal{E}}_{n}(f,f^{q-1})\le t}\frac{1}{n}\Ent_{p,q}(f).
\end{align*}
We now introduce a generalized version of $\Xi_{p,q}^{(n),-1}$. 
\begin{defn}[Generalized R\'enyi-Sobolev Function]
 For $0\le p\le q$ and a function $\beta:[p,q]\to[0,\infty]$, define
a generalized version of (the ``inverse'' of) the R\'enyi-Sobolev
function as 
\begin{equation}
\Upsilon_{p,q}^{(n)}(\beta):=\sup_{f:\frac{1}{(s-1)n}\overline{\mathcal{E}}_{n}(f,f^{s-1})\le\beta(s),\,\forall s\in[p,q]}\frac{1}{n}\Ent_{p,q}(f).\label{eq:Upsilon}
\end{equation}
\end{defn}
Obviously, by choosing $\beta(s)=\begin{cases}
t, & s=q\\
\infty, & s\in[p,q)
\end{cases}$, the quantity $\Upsilon_{p,q}^{(n)}(\beta)$ reduces to $\Xi_{p,q}^{(n),-1}(t)$
with $p\le q$, and by choosing $\beta(s)=\begin{cases}
t, & s=p\\
\infty, & s\in(p,q]
\end{cases}$, the quantity $\Upsilon_{p,q}^{(n)}(\beta)$ reduces to $\Xi_{q,p}^{(n),-1}(t)$
with $p\le q$. Hence, Theorem \ref{thm:RSI} provides  bounds on
$\Upsilon_{p,q}^{(n)}(\beta)$ for these two specific choices of $\beta$.
We next provide an upper bound on $\Upsilon_{p,q}^{(n)}(\beta)$
for a general choice of $\beta$, which is tighter than the bound
in Theorem \ref{thm:RSI} when the constraint in \eqref{eq:Upsilon}
is active for all $s\in[p,q]$.
\begin{thm}[R\'enyi--Sobolev Inequalities]
\label{thm:OU-Renyi-Sobolev} For  any function $\beta:[p,q]\to[0,\infty]$,
\begin{equation}
\Upsilon_{p,q}^{(n)}(\beta)\le\frac{qp}{(q-p)}\int_{p}^{q}\phi_{s}(\beta(s))s^{-2}\d s,\label{eq:-31}
\end{equation}
where  $\phi_{s}:=\conc\,\gamma_{s}$ with 
\begin{equation}
\gamma_{s}(t):=\sup_{f:\frac{1}{s-1}\overline{\mathcal{E}}(f,f^{s-1})\le t}\overline{\Ent}(f^{s})\label{eq:-32}
\end{equation}
and $\conc$ denoting the upper concave envelope is the ``inverse''
of the function $\conv\,\Xi_{s}$ (recall the definition of $\Xi_{s}$
in \eqref{eq:FIpLSI-1}).  The inequality in \eqref{eq:-31} is asymptotically
tight if $\phi_{s}(\beta(s))$ is attained by the same convex combination
of two functions $f_{1},f_{2}$ for all $s\in[p,q]$. 
\end{thm}
\begin{IEEEproof}
  By definition, 
\begin{align*}
\Ent_{p,q}(f) & =\frac{qp}{q-p}\ln\frac{\|f\|_{q}}{\|f\|_{p}},
\end{align*}
which implies that 
\begin{align*}
\frac{\partial}{\partial q}(\frac{1}{p}-\frac{1}{q})\Ent_{p,q}(f) & =\frac{\partial}{\partial q}\ln\frac{\|f\|_{q}}{\|f\|_{p}}=\frac{1}{q^{2}}\overline{\Ent}(f^{q}).
\end{align*}
Hence, 
\begin{align*}
(\frac{1}{p}-\frac{1}{q})\Ent_{p,q}(f) & =\int_{p}^{q}\frac{1}{s^{2}}\overline{\Ent}(f^{s})\d s.
\end{align*}
By the quasi-tensorization property of the log-Sobolev function $\Xi_{q}^{(n)}$
in Theorem \ref{thm:strongSLI}, it holds that under the condition
$\frac{1}{(s-1)n}\overline{\mathcal{E}}_{n}(f,f^{s-1})\le\beta(s)$,
\begin{equation}
\frac{1}{n}\overline{\Ent}(f^{s})\le\phi_{s}(\beta(s)).\label{eq:-33}
\end{equation}
Therefore, 
\begin{align*}
(\frac{1}{p}-\frac{1}{q})\Ent_{p,q}(f) & \le\int_{p}^{q}\frac{n}{s^{2}}\phi_{s}(\beta(s))\d s,
\end{align*}
which yields the desired inequality in \eqref{eq:-31}. 

If there are two functions $f_{1},f_{2}:\calX\to[0,\infty)$ and $\lambda\in[0,1]$
such that for all $s\in[p,q]$, 
\begin{align*}
\frac{1}{s-1}\left(\lambda\overline{\mathcal{E}}(f_{1},f_{1}^{s-1})+(1-\lambda)\overline{\mathcal{E}}(f_{2},f_{2}^{s-1})\right) & =\beta(s),\\
\lambda\overline{\Ent}(f_{1}^{s})+(1-\lambda)\overline{\Ent}(f_{2}^{s}) & =\phi_{s}(\beta(s)),
\end{align*}
then the two sides of \eqref{eq:-33} will asymptotically coincide
for the function $f(x^{n}):=\prod_{i=1}^{\left\lfloor n\lambda\right\rfloor }f_{1}(x_{i})\prod_{i=\left\lfloor n\lambda\right\rfloor +1}^{n}f_{2}(x_{i})$
for any given $s$. This indicates that the inequality in \eqref{eq:-31}
is asymptotically tight. 
\end{IEEEproof}

\subsection{Concentration Inequalities}

We now use Theorem \ref{thm:OU-Renyi-Sobolev} to strengthen existing
concentration inequalities. To this end, we first provide a general
concentration inequality. 
\begin{thm}[General Concentration Inequality]
\label{thm:concentration}Consider a product semigroup with the carr\'e
du champ operator $\Gamma^{\oplus n}$ (defined in \eqref{eq:cdc})
 and stationary distribution $\pi^{n}$. Let $X^{n}\sim\pi^{n}$
and $0\le p\le q$. Let $\beta:[p,q]\to[0,\infty]$ be a given function.
Then, for every function $f:\calX^{n}\to(0,\infty)$ such that 
\begin{equation}
\frac{1}{(s-1)n}\Gamma^{\oplus n}(f,f^{s-1})\le\beta(s)f^{s},\;\forall s\in[p,q],\;\pi^{n}\textrm{-a.s.},\label{eq:-36}
\end{equation}
it holds that 
\begin{align}
 & \mathbb{P}\{\ln f(X^{n})-\ln\|f\|_{p}\ge r\}\nonumber \\
\le & \exp\left\{ nq\int_{p}^{q}\phi_{s}(\beta(s))s^{-2}\d s-rq\right\} ,\;\forall r\in\mathbb{R}.\label{eq:-28}
\end{align}
Here, as usual, $\|f\|_{0}:=\exp\left\{ \mathbb{E}[\ln f(X^{n})]\right\} $. 
\end{thm}
For $p=0$, the probability in \eqref{eq:-28} turns into $\mathbb{P}\{g(X^{n})-\mathbb{E}[g(X^{n})]\ge r\}$,
where $g:=\ln f$. This probability is the common one in the study
of concentration of measure \cite{ledoux2001concentration}.  Moreover,
in the standard setting with $\calX^{n}=\mathbb{R}^{n}$, the condition
$\|\nabla(\ln f)\|_{\ell_{2}}\le1$, $\pi^{n}$-almost surely is usually
assumed, instead of the one in \eqref{eq:-36}. Later we will apply
Theorem \ref{thm:concentration} to the Ornstein--Uhlenbeck and Bonami--Beckner
semigroups, for which we can see that the condition in \eqref{eq:-36}
is a generalization of the standard one.  

Theorem \ref{thm:concentration} holds for arbitrary semigroups, and
hence is a very general result. The proof of Theorem \ref{thm:concentration}
given below resembles the Herbst argument; see e.g., \cite{ledoux2001concentration,RagSason}.
In the standard Herbst argument, standard log-Sobolev inequalities
are used to estimate the $(p,q)$-entropy involved in the Herbst
argument. In contrast, the bound in Theorem \ref{thm:concentration}
(precisely, the bound in Theorem \ref{thm:OU-Renyi-Sobolev}) is derived
via the generalized R\'enyi-Sobolev function $\Upsilon_{p,q}^{(n)}$,
which is estimated by using the nonlinear strengthening of log-Sobolev
inequalities in Theorem \ref{thm:strongSLI}. Hence, it is not surprising
that the bound in Theorem \ref{thm:concentration} is in general tighter
than the ones obtained by the standard Herbst argument, e.g., the
Gaussian concentration bound given in \cite[Section 3.1.4]{RagSason}. 
\begin{IEEEproof}[Proof of Theorem \ref{thm:concentration}]
We only focus on the case $p>0$. By using  the continuity of probability
measures, Theorem \ref{thm:concentration} for $p=0$ follows by letting
$p\downarrow0$.

The normalized Dirichlet form for $f$ satisfies 
\begin{align}
\frac{1}{(s-1)n}\overline{\mathcal{E}}_{n}(f,f^{s-1}) & =\frac{\int\Gamma^{\oplus n}(f,f^{s-1})\d\pi^{n}(x^{n})}{n(s-1)\int f^{s}(x^{n})\d\pi^{n}(x^{n})}\nonumber \\
 & \le\beta(s).
\end{align}
By the Chernoff bound, we obtain that for any $q\ge p$, 
\begin{align}
 & \mathbb{P}\{\ln f-\ln\|f\|_{p}\ge r\}\nonumber \\
 & \le\frac{\mathbb{E}[f^{q}]}{\|f\|_{p}^{q}}\rme^{-rq}\\
 & =\rme^{(q/p-1)\Ent_{p,q}(f)-rq}\label{eq:-22}\\
 & \le\rme^{n(q/p-1)\Upsilon_{p,q}^{(n)}(\beta)-rq}\label{eq:-34}\\
 & \le\rme^{n(q/p-1)\frac{qp}{(q-p)}\int_{p}^{q}\phi_{s}(\beta(s))s^{-2}\d s-rq}\label{eq:-35}\\
 & =\rme^{nq\int_{p}^{q}\phi_{s}(\beta(s))s^{-2}\d s-rq},
\end{align}
where \eqref{eq:-34} follows by definition of $\Upsilon_{p,q}^{(n)}$
in \eqref{eq:Upsilon}, and \eqref{eq:-35} follows by Theorem \ref{thm:OU-Renyi-Sobolev}. 
\end{IEEEproof}

We next apply Theorem \ref{thm:concentration} to recover the Gaussian
concentration inequalities for Gaussian measures.  Denote $\|x^{n}\|:=\|x^{n}\|_{\ell_{2}}=\sqrt{\sum_{i=1}^{n}x_{i}^{2}}$.
Then, $\|x^{n}-y^{n}\|$ is the Euclidean distance between $x^{n},y^{n}$.
We then obtain the following corollary, whose proof is provided in
Section \ref{sec:Proof-of-Corollary}.
\begin{cor}[Concentration for Gaussian Measures]
\label{cor:Gaussian}Let $X^{n}\sim\gamma_{n}:=\mathcal{N}(\boldsymbol{0},\mathbf{I}_{n})$
with $\mathbf{I}_{n}$ denoting the identity matrix of size $n\times n$.
Then, for every differentiable function $g:\mathbb{R}^{n}\to\mathbb{R}$
such that $\|\nabla g(X^{n})\|\le1$ almost surely, it holds that
for $p\ge0$, 
\begin{align}
\mathbb{P}\{g(X^{n})-\ln\|\rme^{g}\|_{p}\ge r\} & \le\begin{cases}
\rme^{-\frac{1}{2}(r+\frac{p}{2})^{2}}, & r\ge\frac{p}{2}\\
\rme^{-pr}, & r<\frac{p}{2}
\end{cases}.\label{eq:-37}
\end{align}
 
\end{cor}
\begin{rem}
This result shows that the probability at the left-hand side obeys
a dimension-free Gaussian concentration bound for $r\ge\frac{p}{2}$
and obeys a dimension-free exponential concentration bound for $r<\frac{p}{2}$.
\end{rem}
By Rademacher's theorem, this corollary still holds
if the assumption that $g$ is differentiable almost everywhere and
$\|\nabla g(X^{n})\|\le1$ almost everywhere is replaced with the
assumption that $g$ is $1$-Lipschitz; see details in the proof of
\cite[Theorem 3.2.2]{RagSason}. The latter is the common assumption
in the study of concentration of measure.  Furthermore, it is common
to investigate the probability $\mathbb{P}\{g(X^{n})-\mathbb{E}[g(X^{n})]\ge r\}$
in the field of concentration of measure. This probability is exactly
the one in \eqref{eq:-37} with $p=0$. For $p=0$, the inequality
in \eqref{eq:-37} reduces to the classic Gaussian concentration bound
\cite{ledoux2001concentration,RagSason}. 

We next apply Theorem \ref{thm:concentration} to derive new concentration
inequalities for the uniform measure on the discrete hypercube. We
then obtain the following corollary, whose proof is provided in Section
\ref{sec:Proof-of-Corollary-1}.
\begin{cor}[Concentration for Uniform Measure on Hypercube]
\label{cor:binary} Let $X^{n}\sim\pi^{n}:=\Unif\{0,1\}^{n}$. Then,
for every function $g:\{0,1\}^{n}\to\mathbb{R}$ such that 
\begin{equation}
|\nabla_{i}g(x^{n})|\le1,\;\forall i\in[n],x^{n}\in\{0,1\}^{n},\label{eq:-26}
\end{equation}
where $\nabla_{i}g(x^{n}):=g(x^{n}\oplus e_{i})-g(x^{n})$ with $x^{n}\oplus e_{i}=(x_{1},...,x_{i-1},1-x_{i},x_{i+1},...,x_{n})$,
it holds that for $p\in[0,2]$, 
\begin{align}
 & \mathbb{P}\{g(X^{n})-\ln\|\rme^{g}\|_{p}\ge r\}\nonumber \\
 & \le\exp\left\{ \inf_{q\in[p,2]}nq\int_{p}^{q}\Xi_{s}^{-1}(\beta(s))s^{-2}\d s-rq\right\} ,\;\forall r\in\mathbb{R},\label{eq:-37-1}
\end{align}
where $\Xi_{s}^{-1}$ is the inverse of the function $\Xi_{s}$ given
in \eqref{eq:-9}, and 
\begin{equation}
\beta(s)=\frac{\left(\rme-1\right)\left(\rme^{(s-1)}-1\right)}{2(s-1)}.\label{eq:-10-1}
\end{equation}
 
\end{cor}
\begin{rem}
The condition in \eqref{eq:-26} is equivalent to saying that  $g$
is $1$-Lipschitz under the Hamming metric $|x^{n}-y^{n}|:=|\{i:x_{i}\neq y_{i}\}|$,
i.e., 
\[
|g(x^{n})-g(y^{n})|\le|x^{n}-y^{n}|,\forall x^{n},y^{n}.
\]
 
\end{rem}
\begin{rem}
It is difficult to obtain a closed-form expression for the bound in
\eqref{eq:-37-1}, but this bound can be evaluated numerically. 
\end{rem}
\begin{rem}
\label{rem:It-can-be}It can be verified that $\beta$ in \eqref{eq:-10-1}
satisfies $\beta(s)\le2$ for all $s\in[0,2]$. So, the bound in \eqref{eq:-37-1}
still holds if $\beta(s)$ is replaced by the constant $2$. 
\end{rem}
The concentration inequality in \eqref{eq:-37-1} can recover the
standard concentration inequality derived by the standard log-Sobolev
inequality. This can be seen as follows. Observe that $\Xi_{s}^{-1}$
is concave and its derivative at $\alpha=0$ is $\frac{s^{2}}{2}$.
So, 
\begin{equation}
\Xi_{s}^{-1}(t)\le\frac{s^{2}t}{2}.\label{eq:-23}
\end{equation}
This is exactly the standard log-Sobolev inequality. Substituting
this inequality into \eqref{eq:-37-1} and setting $\beta(s)$ to
$2$ for simplicity (see Remark  \ref{rem:It-can-be}), it yields
that 
\begin{align}
\mathbb{P}\{g(X^{n})-\ln\|\rme^{g}\|_{p}\ge r\} & \le\rme^{\inf_{q\in[p,2]}nq(q-p)-rq}.\label{eq:-37-1-2}
\end{align}
Since $g$ is $1$-Lipschitz under the Hamming metric, we have that
$|g(X^{n})-\ln\|\rme^{g}\|_{p}|\le n$ a.s. So, it suffices to consider
$r\le n$. For this case, we choose $q$ as the optimal one $q^{*}=\frac{p}{2}+\frac{r}{2n}$
for $r\ge np$, and $q^{*}=p$ for $r<np$, which yields 
\begin{align}
\mathbb{P}\{g(X^{n})-\ln\|\rme^{g}\|_{p}\ge r\} & \le\begin{cases}
\rme^{-n(\frac{r}{2n}+\frac{p}{2})^{2}}, & r\ge np\\
\rme^{-pr}, & r<np
\end{cases}.\label{eq:-24}
\end{align}
In particular, for $p=0$, 
\begin{align}
\mathbb{P}\{g(X^{n})-\mathbb{E}[g(X^{n})]\ge r\} & \le\rme^{-\frac{r^{2}}{4n}},\quad\forall r\ge0,\label{eq:-38}
\end{align}
which is a Gaussian concentration bound.   

Notice that by the strict concavity of $\Xi_{s}^{-1}$, the inequality
in \eqref{eq:-23} is strict for any $t>0$. Hence, the bound in \eqref{eq:-37-1}
is strictly stronger than the one in \eqref{eq:-24}.  In fact,
the bound in \eqref{eq:-37-1} only can improve the factor in the
exponent. This can be seen as follows. Since the function $s\mapsto\Xi_{s}^{-1}(\beta(s))s^{-2}$
is continuous on $[0,2]$, it is also bounded on $[0,2]$. So, the
bound in \eqref{eq:-37-1} is sandwiched between $\rme^{\inf_{q\in[p,2]}C_{1}nq(q-p)-rq}$
and $\rme^{\inf_{q\in[p,2]}C_{2}nq(q-p)-rq}$ for some $0<C_{1}<C_{2}<2$,
which are respectively equal to 
\[
\begin{cases}
\rme^{-C_{1}n(\frac{r}{2C_{1}n}+\frac{p}{2})^{2}}, & r\ge C_{1}np\\
\rme^{-pr}, & r<C_{1}np
\end{cases}
\]
and 
\[
\begin{cases}
\rme^{-C_{2}n(\frac{r}{2C_{2}n}+\frac{p}{2})^{2}}, & r\ge C_{2}np\\
\rme^{-pr}, & r<C_{2}np
\end{cases}.
\]
These expressions look similar to the right-hand side of \eqref{eq:-24},
but with different factors in the exponents. 

Furthermore, it is well known that concentration inequalities can
be rephrased as the estimation of the probability of the $r$-neighborhood
of a set with a given probability \cite{ledoux2001concentration}.
Although here we recover or strengthen the existing concentration
inequalities derived by the standard Herbst argument, there already
exist stronger concentration results in the literature which are derived
by using other methods. For example, for Gaussian measures with the
Euclidean distance, by the symmetrization technique, half-spaces are
proven to be extremizers in the problem of concentration of measure
\cite{sudakov1978extremal,borell1975brunn}; for the uniform measure
with the Hamming distance on the hypercube, by a compression technique,
Hamming balls are proven to be extremizers \cite{bollobas1986combinatorics}.
For other product probability measures, the extremizers are still
unknown in general. However, the exact convergence rate of the target
probability in concentration of measure was already known; see \cite{alon1998asymptotic}
for large sets in finite spaces, \cite{gozlan2005principe,gozlan2007large}
for large sets in arbitrary Polish spaces, and \cite{yu2024exact}
for small sets in arbitrary Polish spaces. The proof in \cite{alon1998asymptotic}
is based on graph theory, the proof in \cite{gozlan2007large} is
based on optimal transport and duality theory, and the proofs in \cite{gozlan2005principe,yu2024exact}
are based on Marton's argument with nonlinear strengthenings of transport-cost
inequalities.  

\section{\protect\label{sec:Connections-to-Spectral}Connections to Spectral
Graph Theory}

We now connect the R\'enyi--Sobolev inequalities to the spectral graph
theory. To this end, we first introduce several new concepts on matrices. 

\subsection{Rayleigh $q$-quotient and Numerical $q$-radius}

We now introduce two new quantities---Rayleigh $q$-quotient and
numerical $q$-radius, which are generalizations of the well-known
Rayleigh quotient and the numerical radius (or the spectral radius). 
\begin{defn}
For a symmetric matrix $\bA\in\mathbb{R}^{\mathcal{X}\times\mathcal{X}}$
and a nonnegative function $f:\mathcal{X}\to[0,\infty)$ that is not
the zero function, we define the Rayleigh $q$-quotient with respect
to $(\bA,f)$ as for $q\in(0,1)\cup(1,\infty)$, 
\[
R_{q}(\bA,f)=\frac{\sum_{x,y}f(x)\bA_{x,y}f(y)^{q-1}}{\sum_{x}f(x)^{q}}.
\]
The Rayleigh $q$-quotient for $q\in\{0,1,\infty\}$ is defined by
the continuous extension. Specifically, for $q\in\{0,1,2,\infty\}$,
\begin{equation}
R_{q}(\bA,f)=\begin{cases}
{\displaystyle \frac{\underset{(x,y)\in\supp(f)^{2}}{\sum}f(x)\bA_{x,y}f(y)^{-1}}{|\supp(f)|}}, & q=0\;\&\\
 & \theta=0,\\
+\infty, & q=0\;\&\\
 & \theta>0,\\
-\infty, & q=0\;\&\\
 & \theta<0,\\
{\displaystyle \frac{\underset{(x,y)\in\supp(f)^{2}}{\sum}f(x)\bA_{x,y}}{\underset{x\in\supp(f)}{\sum}f(x)}}, & q=1,\\
{\displaystyle \frac{\sum_{x,y}f(x)\bA_{x,y}f(y)}{\sum_{x}f(x)^{2}}}, & q=2,\\
{\displaystyle \frac{1}{|\mathcal{X}_{0}|}\sum_{y\in\mathcal{X}_{0}}\sum_{x\in\mathcal{X}}\frac{f(x)}{f_{0}}\bA_{x,y}}, & q=\infty.
\end{cases}\label{eq:-16}
\end{equation}
where $f_{0}:=\max_{x}f(x)$, $\mathcal{X}_{0}:=\{x:f(x)=f_{0}\}$,
and 
\[
\theta:=\sum_{x\in\supp(f),y\notin\supp(f)}f(x)\bA_{x,y}.
\]
\end{defn}
From the expression in \eqref{eq:-16}, it can be seen that when $q=2$,
the Rayleigh $2$-quotient reduces to the standard Rayleigh quotient. 
\begin{defn}
For a nonnegative symmetric matrix $\bA$ and $q\ge0$, define the
(numerical) $q$-radius of $\bA$ as 
\begin{align*}
\rho_{q}(\bA) & =\max_{f}R_{q}(\bA,f),
\end{align*}
where the supremum is taken over all nonnegative functions but not
the zero function. 
\end{defn}
In particular, 
\[
\rho_{q}(\bA)=\begin{cases}
+\infty, & q=0\\
{\displaystyle \max_{x}\sum_{y\in\mathcal{X}}\bA_{x,y}}, & q=1\\
{\displaystyle \max_{f}\frac{\sum_{x,y}f(x)\bA_{x,y}f(y)}{\sum_{x}f(x)^{2}}}, & q=2\\
{\displaystyle \max_{y}\sum_{x\in\mathcal{X}}\bA_{x,y}}, & q=\infty
\end{cases}.
\]
When $q=2$, the numerical $2$-radius is exactly the standard numerical
radius.  By the equivalence between the standard numerical radius
and the spectral radius, the numerical $2$-radius is identical to
the spectral radius. That is, 
\[
\rho_{2}(\bA)=\lambda_{1}(\bA),
\]
where $\lambda_{1}(\bA)$ denotes the maximum eigenvalue of $\bA$.
 Furthermore,  the numerical $2$-radius admits the following nice
properties. 
\begin{prop}
\label{prop:symmetry} For a nonnegative symmetric matrix $\bA$ and
$q\in[1,\infty]$, it holds that 
\begin{equation}
\rho_{q}(\bA)=\rho_{q'}(\bA),\label{eq:-29}
\end{equation}
where $q':=\frac{q}{q-1}$ is the H\"older conjugate of $q$ (with convention
$1'=\infty$ and $\infty'=1$).  Moreover, $\rho_{q}(\bA)$ is nonincreasing
in $q\in[1,2]$ and nondecreasing in $q\in[2,\infty]$. 
\end{prop}
\begin{IEEEproof}
Denote $Q(x):=\frac{f(x)^{q}}{\sum_{x}f(x)^{q}}$. Then, $Q$ is a
distribution on $\mathcal{X}$. Using $Q$, we can rewrite 
\[
R_{q}(\bA,f)=\sum_{x,y}Q^{1/q}(x)\bA_{x,y}Q(y)^{1/q'}.
\]
By symmetry of $\bA$, it is easy to obtain \eqref{eq:-29}. 

Computing the derivative of $R_{q}(\bA,f)$ with respect to $1/q$,
we obtain 
\begin{align*}
 & \frac{\partial}{\partial(1/q)}R_{q}(\bA,f)\\
 & =\sum_{x,y}Q^{1/q}(x)\bA_{x,y}Q(y)^{1/q'}\ln\frac{Q(x)}{Q(y)}\\
 & =\sum_{x<y}\bigl(Q^{1/q}(x)\bA_{x,y}Q(y)^{1/q'}\\
 & \qquad-Q^{1/q}(y)\bA_{x,y}Q(x)^{1/q'}\bigr)\ln\frac{Q(x)}{Q(y)}\\
 & =\sum_{x<y}Q^{1/q}(x)\bA_{x,y}Q(y)^{1/q}\\
 & \qquad\times\left(Q(y)^{1-2/q}-Q(x)^{1-2/q}\right)\ln\frac{Q(x)}{Q(y)}.
\end{align*}
For $q\in[2,\infty]$, it holds that $\ln\frac{Q(x)}{Q(y)}\ge0$ $\Longleftrightarrow$
$Q(x)\ge Q(y)$ $\Longleftrightarrow$ $Q(x)^{1-2/q}\ge Q(y)^{1-2/q}$.
So, $\frac{\partial}{\partial(1/q)}R_{q}(\bA,f)\le0$, which implies
that $R_{q}(\bA,f)$ is nonincreasing in $1/q$ for $q\in[2,\infty]$,
i.e., nondecreasing in $q\in[2,\infty]$. 
\end{IEEEproof}
The $q$-radius also admits an intimate relation with the matrix norm:
\[
\rho_{q}(\bA)=\begin{cases}
\|\bA\|_{1}, & q=1\\
\|\bA\|_{2}, & q=2\\
\|\bA\|_{\infty}, & q=\infty
\end{cases}.
\]
It is worth noting that although the $q$-radius and the matrix $q$-norm
are equal for certain values of $q$, but they are not identical in
general. 

\subsection{Spectral Graph Theory}

We now turn back to connect the R\'enyi--Sobolev inequalities to the
spectral graph theory. Consider a $d$-regular (undirected) graph
$G=(V,E)$. Here $\calX:=V$. We denote $x\sim y$ if $(x,y)$ is
an edge in $G$, i.e., $(x,y)\in E$. Denote $\bA$ as the adjacent
matrix, i.e., $\bA_{x,y}=1$ iff $x\sim y$ in $G$. Obviously, $\bA$
is symmetric. Since $G$ is $d$-regular, $\sum_{y}A_{x,y}=d$. Let
$\bL=\bA-d\bI$, where $\bI$ is the identity matrix. In fact, $-\bL$
is the Laplacian matrix of $G$. The matrix $\bL$ here is an instance
of the general one described at the beginning of the paper. The stationary
distribution $\pi$ for this case is the uniform distribution on $V$. 

Denote $G^{n}$ as the $n$-fold Cartesian graph product of $G$
with itself. In other words, $x^{n}\sim y^{n}$ in $G^{n}$ iff there
is some $i$ such that $x_{i}\sim y_{i}$ in $G$ and $x_{j}=y_{j}$
for all $j\neq i$. Obviously, $G^{n}$ is $nd$-regular and the stationary
distribution with respect to $G^{n}$ is $\pi^{n}$. The normalized
Dirichlet form for this case is 
\begin{align}
\overline{\mathcal{E}}_{n}(f,g) & =\frac{-\langle\Delta f,g\rangle}{\langle f,g\rangle}\\
 & =nd-\frac{\sum_{(x^{n},y^{n}):x^{n}\sim y^{n}}f(x^{n})g(y^{n})}{\sum_{x^{n}}f(x^{n})g(x^{n})},
\end{align}
where $\Delta f(x^{n}):=\sum_{y^{n}:y^{n}\sim x^{n}}(f(y^{n})-f(x^{n}))$.
In particular, 
\begin{align}
\overline{\mathcal{E}}_{n}\big(f,f^{q-1}\big) & =nd-\frac{\sum_{(x^{n},y^{n}):x^{n}\sim y^{n}}f(x^{n})f(y^{n})^{q-1}}{\sum_{x^{n}}f(x^{n})^{q}}\label{eq:Eff-1-1}\\
 & =nd-R_{q}(\bA_{n},f),\label{eq:-17}
\end{align}
where $\bA_{n}$ is the adjacent matrix of $G^{n}$, i.e., the $n$-fold
Kronecker sum of $\bA$ with itself. 

Let $G'=(V',E')$ be an induced subgraph of $G^{n}$. Then, the
adjacent matrix of $G'$ is $\bA_{V',V'}$, which is a submatrix of
$\bA$ with indice restricted to the set $V'\times V'$. 
\begin{defn}
Given $q\ge0$, we define the $q$-radius $\rho_{q}(G')$ of the subgraph
$G'$ as 
\[
\rho_{q}(G'):=\rho_{q}(\bA_{V',V'}).
\]
\end{defn}
By the definition of $\rho_{q}(\bA_{V',V'})$, we have a more explicit
expression for $\rho_{q}(G')$:
\begin{align}
\rho_{q}(G') & =\max_{f}\frac{\sum_{(x^{n},y^{n})\in V'{}^{2}:x^{n}\sim y^{n}}f(x^{n})f(y^{n})^{q-1}}{\sum_{x^{n}\in V'}f(x^{n})^{q}}\\
 & =\max_{f:\supp(f)\subseteq V'}\frac{\sum_{(x^{n},y^{n}):x^{n}\sim y^{n}}f(x^{n})f(y^{n})^{q-1}}{\sum_{x^{n}}f(x^{n})^{q}}\\
 & =\max_{f:\supp(f)\subseteq V'}R_{q}(\bA_{n},f).\label{eq:-16-1}
\end{align}
Here, the maximizations are taken over nonnegative functions but not
the zero function. Proposition \ref{prop:symmetry} immediately implies
$\rho_{q}(G')=\rho_{q'}(G')$, where $q':=\frac{q}{q-1}$ is the H\"older
conjugate of $q$. When $q=2$,  the $2$-radius $\rho_{2}(G')$
of $G'$ is also known as the spectral radius of $G'$, which is equal
to the spectral radius $\lambda_{1}(\bA_{V',V'})$ of the adjacent
matrix. When $q=1$ or $\infty$, $\rho_{1}(G')=\rho_{\infty}(G')$
is equal to the maximum degree of the subgraph $G'$. Furthermore,
combining \eqref{eq:-17} and \eqref{eq:-16-1}, the $q$-radius of
a subgraph and the normalized Dirichlet form admits the following
relation:
\begin{align}
\rho_{q}(G') & =nd-\max_{f:\supp(f)\subseteq V'}\overline{\mathcal{E}}_{n}\big(f,f^{q-1}\big).\label{eq:-18}
\end{align}
When $q=2$, the maximization term in \eqref{eq:-18} is also called
the fractional edge boundary size of $G'$ \cite{samorodnitsky2008modified_article}. 

In spectral graph theory, it has attracted a lot of interest to find
the maximum possible value of the spectral radius of $G'$ over all
subgraphs $G'$ with a given size $|V'|=m$; see, e.g., \cite{friedman2005generalized,bollobas2018eigenvalues,liu2023unsolved}.
This problem is known as the Faber--Krahn problem \cite{friedman2005generalized}.
Here, we generalize this question to the $q$-radius case.
\begin{defn}
Given $q\ge0$ and $1\le m\le|V|^{n}$, the Faber--Krahn maximum
of order $q$ and size $m$ is defined by 
\begin{equation}
\Lambda_{q}(m):=\Lambda_{q}^{(n)}(m):=\max_{G':|V'|\le m}\rho_{q}(G').\label{eq:-18-1}
\end{equation}
\end{defn}
Note that for $q=2$, the Faber--Krahn maximum of order $2$ reduces
to the standard Faber--Krahn maximum. Proposition \ref{prop:symmetry}
immediately implies the following property of $\Lambda_{q}(m)$.
\begin{prop}
\label{prop:symmetry-1} For $n\ge1,q\in[1,\infty]$\textup{,} and
$1\le m\le|V|^{n}$, it holds that 
\[
\Lambda_{q}(m)=\Lambda_{q'}(m),
\]
where $q':=\frac{q}{q-1}$ is the H\"older conjugate of $q$. Moreover,
$\Lambda_{q}(m)$ is nonincreasing in $q\in[1,2]$ and nondecreasing
in $q\in[2,\infty]$. 
\end{prop}
Substituting \eqref{eq:-18} into \eqref{eq:-18-1} yields that 
\begin{align}
\Lambda_{q}(m) & =nd-\min_{G':|V'|\le m}\min_{f:\supp(f)\subseteq V'}\overline{\mathcal{E}}_{n}(f,f^{q-1})\nonumber \\
 & =nd-\min_{f:|\supp(f)|\le m}\overline{\mathcal{E}}_{n}(f,f^{q-1})\nonumber \\
 & =n\left(d-\hat{\Xi}_{0,q}^{(n)}(\alpha)\right),\label{eq:-30}
\end{align}
where $\alpha=\ln|V|-\frac{1}{n}\ln m$ and $\hat{\Xi}_{0,q}^{(n)}$
is a scaled R\'enyi-Sobolev function given by 
\begin{align*}
\hat{\Xi}_{0,q}^{(n)}(\alpha) & :=(q-1)\Xi_{0,q}^{(n)}(\alpha)\\
 & =\inf_{f:\frac{1}{n}\Ent_{0,q}(f)\ge\alpha}\frac{1}{n}\overline{\mathcal{E}}_{n}(f,f^{q-1}).
\end{align*}
Equation \eqref{eq:-30} follows by the definition of the R\'enyi divergence
of order zero and the definition of $\hat{\Xi}_{0,q}^{(n)}$.  Combining
this relation and Theorem \ref{thm:RSI} yields the asymptotically
tight bound for the Faber--Krahn maximum of order $q$. Note that
the asymptotic tightness in the following theorem is our main contribution. 
\begin{thm}
\label{thm:graph}For $q>1,n\ge1$\textup{,} and $1\le m\le|V|^{n}$,
it holds that 
\begin{align}
\frac{1}{n}\Lambda_{q}(m) & =d-\hat{\Xi}_{0,q}^{(n)}(\alpha)\le d-\conv\,\hat{\Xi}_{q}(\alpha),\label{eq:-13}
\end{align}
where $\alpha=\ln|V|-\frac{1}{n}\ln m$ and
\begin{equation}
\hat{\Xi}_{q}(\alpha):=(q-1)\Xi_{q}(\alpha)=\inf_{f:\overline{\Ent}(f^{q})\ge\alpha}\overline{\mathcal{E}}(f,f^{q-1}).\label{eq:-19}
\end{equation}
Moreover, the bound in \eqref{eq:-13} is asymptotically tight for
exponentially large $m$, i.e., $m=e^{n\beta+o(n)}$ with $0\le\beta<\ln|V|$
fixed. In fact, there is a sequence of subgraphs $G_{n}'$ with vertex
set having size no larger than $m$ and formed by the Cartesian product
of at most two $\epsilon$-typical sets, such that  $\lim_{n\to\infty}\frac{1}{n}\rho_{q}(G_{n}')$
could be arbitrarily close to  the upper bound $d-\conv\,\hat{\Xi}_{q}(\ln|V|-\beta)$
as $\epsilon\downarrow0$.  
\end{thm}
\begin{rem}
By Carath\'eodory's theorem, 
\begin{equation}
\conv\,\hat{\Xi}_{q}(\alpha)=\min_{\substack{\alpha_{0},\alpha_{1}\ge0,\lambda\in[0,1]:\\
\lambda\alpha_{0}+(1-\lambda)\alpha_{1}=\alpha
}
}\lambda\hat{\Xi}_{q}(\alpha_{0})+(1-\lambda)\hat{\Xi}_{q}(\alpha_{1}).\label{eq:-20}
\end{equation}
Denote $(\alpha_{0},\alpha_{1},\lambda)$ as the optimal tuple attaining
the minimization in \eqref{eq:-20}. Denote $\hat{f}_{i},i=0,1$
as the optimal functions attaining $\hat{\Xi}_{q}(\alpha_{0})$ and
$\hat{\Xi}_{q}(\alpha_{1})$ respectively.  Denote $Q_{i},i=0,1$
as the distributions such that $\frac{Q_{i}}{\pi}=\frac{\hat{f}_{i}^{q}}{\mathbb{E}_{\pi}[\hat{f}_{i}^{q}]}$,
where $\pi$ is the uniform distribution on $V$. Let $k=\left\lfloor \lambda n\right\rfloor $.
Based on these notations, the vertex set of the asymptotically optimal
subgraph $G_{n}'$ can be written in the form $V_{0}\times V_{1}$,
where $V_{0}$ consists of sequences $x^{k}$ having the empirical
measure almost equal to $Q_{0}$ and $V_{1}$ consists of sequences
$x^{n-k}$ having the same empirical measure almost equal to $Q_{1}$. 
\end{rem}
\begin{rem}
As mentioned in Remark \ref{rem:The-log-Sobolev-function}, the lower
convex envelope in \eqref{eq:-13} cannot be removed. 
\end{rem}

\begin{IEEEproof}[Proof of Theorem \ref{thm:graph}]
 By Theorem \ref{thm:RSI}, $\hat{\Xi}_{0,q}^{(n)}(\alpha)\ge\conv\,\hat{\Xi}_{q}(\alpha),$
which implies \eqref{eq:-13}. The asymptotic tightness follows by
\eqref{eq:asympt_tight}. 
\end{IEEEproof}
In particular, for the binary case (i.e., $\mathcal{X}=\{0,1\}$ with
$d=1$) and $q>1$, 
\begin{align}
\hat{\Xi}_{q}(\alpha) & =1\!-\!y^{\frac{1}{q}}(1-y)^{\frac{1}{q'}}\!-\!y^{\frac{1}{q'}}(1-y)^{\frac{1}{q}},\label{eq:-9-1-1}
\end{align}
where $y=h^{-1}(\ln2-\alpha)$. Note that compared to the expression
in \eqref{eq:-9}, the factor $1/2$ is removed in the expression
in \eqref{eq:-9-1-1}. This is because the matrix $\bL$ in this subsection
is twice the one in Section \ref{subsec:Binary-Example}. 

The function $\hat{\Xi}_{q}$ for this binary case is convex \cite{polyanskiy2019improved}.
So, 
\begin{align}
\frac{1}{n}\Lambda_{q}(m) & \le1-\,\hat{\Xi}_{q}(\alpha)\nonumber \\
 & =y^{\frac{1}{q}}(1-y)^{\frac{1}{q'}}\!+\!y^{\frac{1}{q'}}(1-y)^{\frac{1}{q}},\label{eq:-21}
\end{align}
where $y=h^{-1}(\ln2-\alpha)=h^{-1}(\frac{1}{n}\ln m)$. This upper
bound is asymptotically tight for exponentially large $m$. A sequence
of subgraphs with vertex sets constituted by Hamming spherical shells
(or Hamming balls) of certain radii asymptotically attains the upper
bound.  When $q=2$, \eqref{eq:-21} reduces to 
\begin{align*}
\frac{1}{n}\Lambda_{2}(m) & \le2\sqrt{y(1-y)},
\end{align*}
which recovers an existing result proven by Samorodnitsky \cite{samorodnitsky2008modified_article}
and Bollob\'as, Lee, and Letzter \cite{bollobas2018eigenvalues}.
The asymptotic tightness of this special case was first proven by
Friedman and Tillich \cite{friedman2005generalized} using a spectral
graph-theoretic method, which is different from the information-theoretic
method used by us. Although the asymptotic optimality is already
known, it is still open to determine the exact value of $\Lambda_{2}(m)$
(or more generally, $\Lambda_{q}(m)$) and the corresponding optimal
sets for every given pair $(n,m)$ \cite{bollobas2018eigenvalues,liu2023unsolved}.
Note that Hamming balls are not the exact solution to this open problem,
at least for certain cases. For example, in the balanced case, i.e.,
$m=2^{n-1}$, $(n-1)$-dimensional subcubes are optimal, but the Hamming
balls of radius $n/2$ are not; see the comments below Theorem 1.4
in \cite{samorodnitsky2008modified_article}. 

As one of important applications, the standard Faber--Krahn maximum
was applied by Friedman and Tillich \cite{friedman2005generalized}
and also by Navon and Samorodnitsky \cite{navon2009linear} to provide
an alternative proof for the famous McEliece--Rodemich--Rumsey--Welch
(MRRW) bound in coding theory via a connection between packing bounds
and isoperimetry. Specifically, they converted the task of finding
an upper bound on the size of a code with the minimum distance given
into the task of finding a lower bound for the Faber--Krahn problem.
Intuitively, it is possible to use our Theorem \ref{thm:graph} generalize
their proofs from the binary case (i.e., the complete graph $K_{2}$)
to arbitrary regular graphs. It might yield some interesting bounds
on the independence number of the product of an arbitrary regular
graph. This is left as a future work. 

\section{\protect\label{sec:Proof-of-Theorem}Proof of Theorem \ref{thm:RSI}}

Before giving the proof of Theorem \ref{thm:RSI}, we first provide
information-theoretic characterizations of $\Xi_{p,q}$ and $\Xi_{p,q}^{(n)}$.
For any $f$, we can choose a distribution $Q$ such that $\frac{f^{q}}{\mathbb{E}[f^{q}]}=\frac{Q}{\pi}$.
From this relationship, we obtain the following information-theoretic
characterizations. Recall that $q'=q/(q-1)$ is the H\"older conjugate
of $q$. 
\begin{prop}
\label{prop:It-holds-that}It holds that for $p\in\mathbb{R}$ and
$q\in\mathbb{R}\backslash\{0\}$, the R\'enyi-Sobolev function satisfies
\begin{align*}
\Xi_{p,q}(\alpha) & =\inf_{Q_{X}:D_{p/q}(Q_{X}\|\pi)\ge\alpha}\frac{1}{q-1}\\
 & \qquad\times\mathcal{E}\left(\Big(\frac{Q_{X}}{\pi}\Big)^{1/q},\Big(\frac{Q_{X}}{\pi}\Big)^{1/q'}\right),
\end{align*}
and in particular, for $p=q$, 
\begin{align}
\Xi_{q}(\alpha) & =\inf_{Q_{X}:D(Q_{X}\|\pi)\ge\alpha}\frac{1}{q-1}\nonumber \\
 & \qquad\times\mathcal{E}\left(\Big(\frac{Q_{X}}{\pi}\Big)^{1/q},\Big(\frac{Q_{X}}{\pi}\Big)^{1/q'}\right).\label{eq:-11}
\end{align}
Similarly, for $p\in\mathbb{R}$ and $q\in\mathbb{R}\backslash\{0\}$,
the $n$-dimensional R\'enyi-Sobolev function satisfies
\begin{align}
\Xi_{p,q}^{(n)}(\alpha) & =\inf_{Q_{X^{n}}:\frac{1}{n}D_{p/q}(Q_{X^{n}}\|\pi^{n})\ge\alpha}\frac{1}{(q-1)n}\nonumber \\
 & \qquad\times\mathcal{E}_{n}\left(\Big(\frac{Q_{X^{n}}}{\pi^{n}}\Big)^{1/q},\Big(\frac{Q_{X^{n}}}{\pi^{n}}\Big)^{1/q'}\right).\label{eq:ITCh}
\end{align}
\end{prop}
\begin{IEEEproof}
If we write ${f^{q}}/{\mathbb{E}[f^{q}]}={Q_{X}}/{\pi}$ for a distribution
$Q_{X}\ll\pi$, then 
\begin{align}
\Ent_{p,q}(f) & =\Ent_{p/q}(f^{q})=D_{p/q}(Q_{X}\|\pi),\label{eq:FI-3-1-1}\\
\overline{\mathcal{E}}(f,f^{q-1}) & =\mathcal{E}\left(\Big(\frac{Q_{X}}{\pi}\Big)^{1/q},\Big(\frac{Q_{X}}{\pi}\Big)^{1/q'}\right).\label{eq:FI-7-1-1}
\end{align}
Uniting \eqref{eq:FI-3-1-1} and \eqref{eq:FI-7-1-1}, one can obtain
the information-theoretic characterization of~ $\Xi_{p,q}$. The
information-theoretic characterization of~ $\Xi_{p,q}^{(n)}$ follows
similarly. 
\end{IEEEproof}
We next prove Theorem \ref{thm:RSI}. 

\subsection{Case of $p\le q$}

We now prove \eqref{eq:RenyiSobolev} and \eqref{eq:asympt_tight}
for $p\le q$. Since given $f$, $\Ent_{p,q}(f)$ is non-decreasing
in $p$ (or $q$) given $q$ (or $p$) \cite{yu2021strong_article},
it holds that $\Ent_{p,q}(f)\le\Ent_{q,q}(f)=\Ent_{q}(f)$ for $p\le q$.
(This result can be alternatively obtained from the monotonicity of
R\'enyi divergence in its order.) Hence, substituting this into the
definition of $\Xi_{p,q}^{(n)}$ yields that  $\Xi_{p,q}^{(n)}(\alpha)\ge\Xi_{q}^{(n)}(\alpha)\ge\conv\,\Xi_{q}(\alpha)$
for $p\le q$, where Theorem \ref{thm:strongSLI} is applied. 

We next use the information-theoretic characterization in Proposition
\ref{prop:It-holds-that} to prove the asymptotic tightness of this
bound for $q\ge1$, i.e., \eqref{eq:asympt_tight}.  Since $\Xi_{p,q}^{(n)}$
is nonincreasing in $p$, it suffices to prove $\Xi_{0,q}^{(n)}(\alpha)\to\conv\,\Xi_{q}(\alpha)$
as $n\to\infty$. 

Let $R$ be an optimal distribution attaining the infimization in
the expression of $\Xi_{q}(\alpha)$ in \eqref{eq:-11}.  Let $T_{x^{n}}$
denote the empirical measure of a sequence $x^{n}$. Let $\mathcal{T}_{\epsilon}^{(n)}:=\mathcal{T}_{\epsilon}^{(n)}(R):=\{x^{n}:|T_{x^{n}}(x)-R(x)|\le\epsilon R(x),\forall x\}$
be the $\epsilon$-typical set of $R$.  We construct a distribution
$Q_{X^{n}}:=R^{n}(\cdot|\mathcal{T}_{\epsilon}^{(n)})$ which is a
conditional distribution of $R^{n}$ on the $\epsilon$-typical set. 

 For a given $Q_{X^{n}}$, define the conditional Radon--Nikodym
derivative of the $k$-th coordinate as 
\begin{equation}
\ell_{k}(y|x^{\setminus k}):=\frac{Q_{X_{k}|X^{\setminus k}}(y|x^{\setminus k})}{\pi(y)},\;\forall\,(y,x^{\setminus k})\in\calX\times\calX^{n-1}.\label{eqn:zeta_k}
\end{equation}
Then, it follows that 
\begin{align}
 & \mathcal{E}_{n}\left(\Big(\frac{Q_{X^{n}}}{\pi^{n}}\Big)^{1/q},\Big(\frac{Q_{X^{n}}}{\pi^{n}}\Big)^{1/q'}\right)\nonumber \\*
 & =-\sum_{k=1}^{n}\sum_{x^{\setminus k}}Q_{X^{\setminus k}}(x^{\setminus k})\nonumber \\
 & \quad\times\left(\sum_{x,y}\ L_{x,y}\ell_{k}(y|x^{\setminus k})^{1/q}\ \ell_{k}(x|x^{\setminus k})^{1/q'}\ \pi(x)\right)\label{eq:-14}\\
 & =\sum_{k=1}^{n}\sum_{x^{\setminus k}}Q_{X^{\setminus k}}(x^{\setminus k})\mathcal{E}\left(\ell_{k}(\cdot|x^{\setminus k})^{1/q},\ell_{k}(\cdot|x^{\setminus k})^{1/q'}\right)\\
 & =\sum_{k=1}^{n}\sum_{x^{\setminus k}}Q_{X^{\setminus k}}(x^{\setminus k})\nonumber \\
 & \quad\times\mathcal{E}\left(\Big(\frac{Q_{X_{k}|X^{\setminus k}=x^{\setminus k}}}{\pi}\Big)^{1/q},\Big(\frac{Q_{X_{k}|X^{\setminus k}=x^{\setminus k}}}{\pi}\Big)^{1/q'}\right),\label{eqn:use_calE}
\end{align}
where \eqref{eq:-14} follows from the definition of the Dirichlet
form in~\eqref{eqn:dirichlet_form}. 

For any sequence $x^{n}$ and any $a\in\mathcal{X}$, 
\begin{align*}
 & |T_{x^{n}}(a)-T_{x^{n-1}}(a)|\\
 & =\left|\frac{1}{n(n-1)}\left(\sum_{i=1}^{n-1}1\{x_{i}=a\}\right)-\frac{1}{n}1\{x_{n}=a\}\right|\\
 & \le\frac{1}{n}.
\end{align*}
By symmetry, for any $k\in[n]$, 
\begin{equation}
|T_{x^{n}}(a)-T_{x^{\backslash k}}(a)|\le\frac{1}{n}.\label{eq:}
\end{equation}

Let $\epsilon_{1}>\epsilon>\epsilon_{0}$ and $\mathcal{T}_{\epsilon_{0}}^{(n-1)}$
be a $\epsilon_{0}$-typical set of length-$(n-1)$ sequences for
$R$. Then, by \eqref{eq:}, for sufficiently large $n$, $x^{\setminus k}\in\mathcal{T}_{\epsilon_{0}}^{(n-1)}$
implies $x^{n}\in\mathcal{T}_{\epsilon}^{(n)}$ for $x_{k}$ taking
any values. Denote $Q_{X^{\setminus k}}$ as the $X^{\setminus k}$-marginal
of $Q_{X^{n}}.$ For $x^{\setminus k}\in\mathcal{T}_{\epsilon_{0}}^{(n-1)}$
and sufficiently large $n$ (which means $x^{n}\in\mathcal{T}_{\epsilon}^{(n)}$),
the marginal probability $Q_{X^{\setminus k}}(x^{\setminus k})$ can
be expressed as follows: 
\begin{align}
Q_{X^{\setminus k}}(x^{\setminus k}) & =\sum_{x_{k}}R^{n}(x^{n}|\mathcal{T}_{\epsilon}^{(n)})\nonumber \\
 & =\sum_{x_{k}}\frac{R^{n}(x^{n})1\{x^{n}\in\mathcal{T}_{\epsilon}^{(n)}\}}{R^{n}(\mathcal{T}_{\epsilon}^{(n)})}\nonumber \\
 & =\sum_{x_{k}}\frac{R^{n}(x^{n})}{R^{n}(\mathcal{T}_{\epsilon}^{(n)})}=\frac{R^{n-1}(x^{\setminus k})}{R^{n}(\mathcal{T}_{\epsilon}^{(n)})}.\label{eq:Q}
\end{align}
Therefore, the conditional probability $Q_{X_{k}|X^{\setminus k}}(x_{k}|x^{\setminus k})$
can be expressed as follows: 
\begin{align}
Q_{X_{k}|X^{\setminus k}}(x_{k}|x^{\setminus k}) & =\frac{Q_{X^{n}}(x^{n})}{Q_{X^{\setminus k}}(x^{\setminus k})}\nonumber \\
 & =\frac{R^{n}(x^{n})}{R^{n}(\mathcal{T}_{\epsilon}^{(n)})}\Big/\frac{R^{n-1}(x^{\setminus k})}{R^{n}(\mathcal{T}_{\epsilon}^{(n)})}=R(x_{k}).\label{eq:-1}
\end{align}
This result implies that if $n$ is sufficiently large and $x^{\setminus k}\in\mathcal{T}_{\epsilon_{0}}^{(n-1)}$,
the conditional distribution $Q_{X_{k}|X^{\setminus k}=x^{\setminus k}}$
is exactly the original distribution $R$ used to define $Q_{X^{n}}$. 

We next need the following fact. 
\begin{fact}
\label{fact:For-,-} For $q\ge1$, $\mathcal{E}\left(\Big(\frac{S}{\pi}\Big)^{1/q},\Big(\frac{S}{\pi}\Big)^{1/q'}\right)$
is uniformly bounded for any distribution $S$. 
\end{fact}
This fact is obvious since it holds that 
\begin{align*}
 & \mathcal{E}\left(\Big(\frac{S}{\pi}\Big)^{1/q},\Big(\frac{S}{\pi}\Big)^{1/q'}\right)\\
 & =-\sum_{x,y}L_{x,y}\Big(\frac{S(y)}{\pi(y)}\Big)^{1/q}\Big(\frac{S(x)}{\pi(x)}\Big)^{1/q'}\pi(x).
\end{align*}
The RHS above is continuous in $S$ on  the probability simplex.
So, the term $S\mapsto\mathcal{E}\left(\frac{S}{\pi},\ln\frac{S}{\pi}\right)$
is bounded. 

By the fact above, there is a number $M$ such that for any $S$,
\begin{equation}
\left|\mathcal{E}\left(\Big(\frac{S}{\pi}\Big)^{1/q},\Big(\frac{S}{\pi}\Big)^{1/q'}\right)\right|\le M.\label{eq:-2}
\end{equation}
Then, substituting \eqref{eq:-1} and \eqref{eq:-2} into \eqref{eqn:use_calE},
we obtain that 
\begin{align}
 & \mathcal{E}_{n}\left(\Big(\frac{Q_{X^{n}}}{\pi^{n}}\Big)^{1/q},\Big(\frac{Q_{X^{n}}}{\pi^{n}}\Big)^{1/q'}\right)\nonumber \\
 & \in\sum_{k=1}^{n}p_{k}\mathcal{E}\left(\Big(\frac{R}{\pi}\Big)^{1/q},\Big(\frac{R}{\pi}\Big)^{1/q'}\right)\nonumber \\
 & \qquad\pm\sum_{k=1}^{n}(1-p_{k})M,\label{eq:-8}
\end{align}
where $p_{k}:=Q_{X^{\setminus k}}(\mathcal{T}_{\epsilon_{0}}^{(n-1)})$. 

We now estimate $p_{k}$. By symmetry, $p_{k}$ are the same for all
$k$. So, we only consider $p_{n}$. For sufficiently large $n$,
\begin{align}
p_{n} & =\sum_{x^{n-1}\in\mathcal{T}_{\epsilon_{0}}^{(n-1)}}\sum_{x_{n}}\frac{R^{n}(x^{n})1\{x^{n}\in\mathcal{T}_{\epsilon}^{(n)}\}}{R^{n}(\mathcal{T}_{\epsilon}^{(n)})}\nonumber \\
 & =\sum_{x^{n-1}\in\mathcal{T}_{\epsilon_{0}}^{(n-1)}}\frac{R^{n-1}(x^{n-1})}{R^{n}(\mathcal{T}_{\epsilon}^{(n)})}\label{eq:-4}\\
 & =\frac{R^{n-1}(\mathcal{T}_{\epsilon_{0}}^{(n-1)})}{R^{n}(\mathcal{T}_{\epsilon}^{(n)})},\nonumber 
\end{align}
where \eqref{eq:-4} follows since $x^{n}\in\mathcal{T}_{\epsilon}^{(n)}$
always holds as long as $x^{n-1}\in\mathcal{T}_{\epsilon_{0}}^{(n-1)}$. 

As an well-known result (called the asymptotic equipartition property)
in information theory, the probability of a typical set converges
to one (exponentially fast) as the dimension goes to infinity. That
is, both $R^{n}(\mathcal{T}_{\epsilon}^{(n)})$ and $R^{n-1}(\mathcal{T}_{\epsilon_{0}}^{(n-1)})$
converge to one. This implies $p_{n}\to1$. Therefore, \eqref{eq:-8}
implies 
\begin{align}
 & \frac{1}{n}\mathcal{E}_{n}\left(\Big(\frac{Q_{X^{n}}}{\pi^{n}}\Big)^{1/q},\Big(\frac{Q_{X^{n}}}{\pi^{n}}\Big)^{1/q'}\right)\nonumber \\
 & \to\mathcal{E}\left(\Big(\frac{R}{\pi}\Big)^{1/q},\Big(\frac{R}{\pi}\Big)^{1/q'}\right).\label{eq:-5}
\end{align}

We next estimate the R\'enyi divergence for the case $p=0$. Observe
that 
\begin{align}
D_{0}(Q_{X^{n}}\|\pi^{n}) & =-\ln\sum_{x^{n}:Q_{X^{n}}(x^{n})>0}\pi^{n}(x^{n})\nonumber \\
 & =-\ln\sum_{x^{n}\in\mathcal{T}_{\epsilon}^{(n)}(R)}e^{n\sum_{x}T_{x^{n}}(x)\log\pi(x)}\nonumber \\
 & \ge-\ln\sum_{x^{n}\in\mathcal{T}_{\epsilon}^{(n)}(R)}e^{n\left(\sum_{x}R(x)\log\pi(x)+o_{\epsilon}(1)\right)}\label{eq:-12}\\
 & \ge-\ln\rme^{n\left(H(R)+\sum_{x}R(x)\log\pi(x)+o_{\epsilon}(1)\right)}\nonumber \\
 & =n\left(D(R\|\pi)+o_{\epsilon}(1)\right),\nonumber 
\end{align}
where $o_{\epsilon}(1)$ denotes a term vanishing as $\epsilon\to0$,
and \eqref{eq:-12} follows by the typical average lemma \cite{elgamal}.

Letting $n\to\infty$ first and $\epsilon\to0$ then, we obtain that
\begin{align}
\frac{1}{n}D_{0}(Q_{X^{n}}\|\pi^{n}) & \to D(R\|\pi).\label{eq:-6}
\end{align}
Combining \eqref{eq:-5}, \eqref{eq:-6}, and the information-theoretic
characterization of $\Xi_{p,q}^{(n)}$ in \eqref{eq:ITCh} yields
that for any $\delta>0$, 
\[
\lim_{n\to\infty}\Xi_{0,q}^{(n)}(\alpha)\le\Xi_{q}(\alpha+\delta).
\]
By the time-sharing argument, we easily obtain that 
\[
\lim_{n\to\infty}\Xi_{0,q}^{(n)}(\alpha)\le\conv\,\Xi_{q}(\alpha+\delta).
\]
Letting $\delta\to0$ and utilizing the continuity of a convex function
on an open interval, we obtain that for $\alpha\in[0,\ln\frac{1}{\min_{x}\pi(x)})$,
\begin{equation}
\lim_{n\to\infty}\Xi_{0,q}^{(n)}(\alpha)\le\conv\,\Xi_{q}(\alpha).\label{eq:-7}
\end{equation}
This, combined with \eqref{eq:RenyiSobolev}, proves \eqref{eq:asympt_tight}
for the case $p=0$, which in turn implies the case $p\le q$.

\subsection{Case of $p>q$}

We next prove \eqref{eq:RenyiSobolev} and \eqref{eq:asympt_tight}
for $p>q$. For this case, from \eqref{eqn:use_calE} and \eqref{eq:positivity},
it is easily seen that 
\[
\frac{1}{(q-1)n}\mathcal{E}_{n}\left(\Big(\frac{Q_{X^{n}}}{\pi^{n}}\Big)^{1/q},\Big(\frac{Q_{X^{n}}}{\pi^{n}}\Big)^{1/q'}\right)\ge0.
\]
That is, 
\[
\Xi_{p,q}^{(n)}(\alpha)\ge0.
\]
We next prove the asymptotic tightness of this bound for $q\ge1$. 

Denote $z$ as $\pi(z)=\min_{x}\pi(x)$, and $z^{n}=(z,z,...,z)$.
   Let $\epsilon>\epsilon_{0}$. Let $\mathcal{T}_{\epsilon}^{(n)}$
and $\mathcal{T}_{\epsilon_{0}}^{(n-1)}$ be the $\epsilon$-typical
set of length-$n$ and length-$(n-1)$ sequences for $\pi$. Since
$\pi\neq\delta_{z}$ where $\delta_{z}$ is the Dirac measure at $z$,
for sufficiently small $\epsilon$, it holds that $z^{n}\notin\mathcal{T}_{\epsilon}^{(n)}$.
We choose $\epsilon$ small enough to satisfy this condition. 

Let $\beta>0$. Denote  $Q_{X^{n}}:=(1-e^{-n\beta})\pi^{n}(\cdot|\mathcal{T}_{\epsilon}^{(n)})+e^{-n\beta}\delta_{z^{n}}$,
where $\delta_{z^{n}}$ is the Dirac measure at $z^{n}$. For sufficiently
large $n$ and $x^{\setminus k}\in\mathcal{T}_{\epsilon_{0}}^{(n-1)}$,
it holds that $x^{n}\in\mathcal{T}_{\epsilon}^{(n)}$; see the arguments
in the proof of the case $p\le q$. We hence have that  $x^{n}\neq z^{n}$
(i.e., $\delta_{z^{n}}(x^{n})=0$) for sufficiently large $n$, since
$z^{n}\notin\mathcal{T}_{\epsilon}^{(n)}$. 

By the fact $\delta_{z^{n}}(x^{n})=0$ and by \eqref{eq:Q}, it holds
that 

\begin{align*}
Q_{X^{\setminus k}}(x^{\setminus k}) & =\sum_{x_{k}}(1-\rme^{-n\beta})\pi^{n}(x^{n}|\mathcal{T}_{\epsilon}^{(n)})\\
 & =(1-\rme^{-n\beta})\frac{\pi^{n-1}(x^{\setminus k})}{\pi^{n}(\mathcal{T}_{\epsilon}^{(n)})}.
\end{align*}
Therefore, 
\begin{align*}
Q_{X_{k}|X^{\setminus k}}(x_{k}|x^{\setminus k}) & =\frac{Q_{X^{n}}(x^{n})}{Q_{X^{\setminus k}}(x^{\setminus k})}\\
 & =\frac{(1-\rme^{-n\beta})\frac{\pi^{n}(x^{n})}{\pi^{n}(\mathcal{T}_{\epsilon}^{(n)})}}{(1-\rme^{-n\beta})\frac{\pi^{n-1}(x^{\setminus k})}{\pi^{n}(\mathcal{T}_{\epsilon}^{(n)})}}\\
 & =\pi(x_{k}).
\end{align*}
Following same proof steps, \eqref{eq:-5} still holds for this case. 

On the other hand, for sufficiently large $n$, 
\begin{align*}
 & D_{p/q}(Q_{X^{n}}\|\pi^{n})\\
 & =\frac{1}{p/q-1}\ln\sum_{x^{n}}Q_{X^{n}}(x^{n})\left(\frac{Q_{X^{n}}(x^{n})}{\pi^{n}(x^{n})}\right)^{p/q-1}\\
 & =\frac{1}{p/q-1}\ln\Biggl(\sum_{x^{n}\in\mathcal{T}_{\epsilon}^{(n)}}\left((1-\rme^{-n\beta})\frac{\pi^{n}(x^{n})}{\pi^{n}(\mathcal{T}_{\epsilon}^{(n)})}\right)^{p/q}\\
 & \qquad\times\left(\frac{1}{\pi^{n}(x^{n})}\right)^{p/q-1}+\rme^{-n\beta p/q}\left(\frac{1}{\pi^{n}(z^{n})}\right)^{p/q-1}\Biggr)\\
 & \ge n\left(\ln\frac{1}{\pi(z)}-\beta\frac{p/q}{p/q-1}\right),
\end{align*}
where in the last line the bound is obtained by discarding the first
term in the logarithm above. By choosing $\beta$ small enough, we
have that the limit (or the limit inferior) of $\frac{1}{n}D_{p/q}(Q_{X^{n}}\|\pi^{n})$
can be greater than or arbitrarily close to $\ln\frac{1}{\pi(z)}$.
Combining this with \eqref{eq:-5} with $R=\pi$ implies that $\lim_{n\to\infty}\Xi_{p,q}^{(n)}(\alpha)\le0$
for all $\alpha\in[0,-\ln\min_{x}\pi(x))$.

\section{\protect\label{sec:Proof-of-Corollary}Proof of Corollary \ref{cor:Gaussian}}

Consider $\mathcal{X}^{n}=\mathbb{R}^{n}$, $\pi^{n}=\gamma_{n}$,
and $T_{t}^{\otimes n}$ is the Ornstein--Uhlenbeck operator given
by 
\begin{align}
 & T_{t}^{\otimes n}f(x^{n})\nonumber \\
 & =\int f(\rme^{-t}x^{n}+\sqrt{1-\rme^{-2t}}z^{n})\d\gamma_{n}(z^{n}),\;\forall x^{n}\in\mathbb{R}^{n}.
\end{align}
For the Ornstein--Uhlenbeck semigroup, it is well known \cite{bakry2004functional}
that the corresponding carr\'e du champ operator is 
\begin{align}
\Gamma^{\oplus n}(f,g) & =\nabla f\cdot\nabla g,\;\forall f,g.\label{eq:FIdiri_cube-1}
\end{align}

Let $g:\mathbb{R}^{n}\to\mathbb{R}$ be such that $\|\nabla g(X^{n})\|\le1$
almost surely. Denote $f=\rme^{g}$. Then, obviously,
\[
\nabla g(x^{n})=\frac{\nabla f(x^{n})}{f(x^{n})}.
\]
So, it holds that
\begin{equation}
\|\nabla f(X^{n})\|\le f(X^{n}),\;\textrm{almost surely}.\label{eq:-25}
\end{equation}
The carr\'e du champ operator satisfies
\begin{align*}
\frac{1}{s-1}\Gamma^{\oplus n}(f,f^{s-1})(X^{n}) & =\|\nabla f(X^{n})\|^{2}f^{s-2}(X^{n})\\
 & \le f^{s}(X^{n}),\;\textrm{a.s.}
\end{align*}
So, $f$ satisfies the assumption in Theorem \ref{thm:concentration}
with $\beta(s)=\frac{1}{n}$, and hence, the inequality in \eqref{eq:-28}
holds. 

We now compute the function $\phi_{s}$ for the Ornstein--Uhlenbeck
semigroup. By the Gaussian log-Sobolev inequality \cite{gross1975logarithmic},
for $s>1$, 
\[
\overline{\Ent}(f^{s})\le\frac{s^{2}}{2(s-1)}\overline{\mathcal{E}}(f,f^{s-1}),\quad\forall f\ge0.
\]
So, $\gamma_{s}(t)\le\frac{s^{2}t}{2}.$ By setting $f$ to be $\rme^{bx-b^{2}/2}$
for a proper $b$, one can see that this upper bound is tight. That
is, $\gamma_{s}(t)=\frac{s^{2}t}{2}$, and hence, 
\[
\phi_{s}(t)=\frac{s^{2}t}{2}.
\]

By Theorem \ref{thm:concentration} with $\beta(s)=\frac{1}{n}$ and
$\phi_{s}(t)=\frac{s^{2}t}{2}$, it holds that 
\begin{align}
\mathbb{P}\{\ln f-\ln\|f\|_{p}\ge r\} & \le\rme^{\frac{q(q-p)}{2}-rq}.\label{eq:-28-1}
\end{align}
Since $q\ge p$ is arbitrary, we choose $q$ as the optimal one $q^{*}=\frac{p}{2}+r$
for $r\ge\frac{p}{2}$ and $q^{*}=p$ for $r<\frac{p}{2}$, which
yields 
\begin{align}
\mathbb{P}\{\ln f-\ln\|f\|_{p}\ge r\} & \le\begin{cases}
\rme^{-\frac{1}{2}(r+\frac{p}{2})^{2}}, & r\ge\frac{p}{2}\\
\rme^{-pr}, & r<\frac{p}{2}
\end{cases}.
\end{align}

\section{\protect\label{sec:Proof-of-Corollary-1}Proof of Corollary \ref{cor:binary}}

Consider  the Bonami--Beckner operator $T_{t}^{\otimes n}$ with
$T_{t}$ given by 
\begin{equation}
T_{t}f(y)=f(y)\frac{1+\rme^{-t}}{2}+f(1-y)\frac{1-\rme^{-t}}{2},\;y\in\{0,1\}.\label{eq:FIhypercube-2}
\end{equation}
For this semigroup, it is well known \cite{bakry2004functional}
that the corresponding carr\'e du champ operator is 
\begin{align*}
 & \Gamma^{\oplus n}(f,g)(x^{n})\\
 & =\frac{1}{2}\sum_{y^{n}\sim x^{n}}(f(y^{n})-f(x^{n}))(g(y^{n})-g(x^{n})),\;\forall f,g,
\end{align*}
where $y^{n}\sim x^{n}$ means that $y^{n},x^{n}\in\{0,1\}^{n}$ differ
in exactly one coordinate. 

Let $g$ be a function satisfying the condition in \eqref{eq:-26}.
Denote $f=\rme^{g}$. Then, for $s\in(0,1)\cup(1,2]$, 
\begin{align*}
 & \frac{1}{s-1}\Gamma^{\oplus n}(f,f^{s-1})(x^{n})\\
 & =\frac{1}{2(s-1)}\sum_{y^{n}\sim x^{n}}(f(y^{n})-f(x^{n}))(f^{s-1}(y^{n})-f^{s-1}(x^{n}))\\
 & =\frac{f^{s}(x^{n})}{2}\sum_{y^{n}\sim x^{n}}\left(\rme^{g(y^{n})-g(x^{n})}-1\right)\\
 & \qquad\times\left(\frac{\rme^{(s-1)[g(y^{n})-g(x^{n})]}-1}{s-1}\right)\\
 & \le\frac{f^{s}(x^{n})}{2}\sum_{y^{n}\sim x^{n}}\left(\rme-1\right)\left(\frac{\rme^{|s-1|}-1}{|s-1|}\right)\\
 & =\frac{\left(\rme-1\right)\left(\rme^{(s-1)}-1\right)}{2(s-1)}nf^{s}(x^{n}).
\end{align*}
By taking limits, this inequality still holds for $s\in\{0,1\}$.
 So, $f$ satisfies the assumption in Theorem \ref{thm:concentration}
with 
\begin{equation}
\beta(s)=\frac{\left(\rme-1\right)\left(\rme^{(s-1)}-1\right)}{2(s-1)}\label{eq:-10}
\end{equation}
and $0\le p\le q\le2$, and hence, the inequality in \eqref{eq:-28}
holds. 

For the Bonami--Beckner semigroup, the function $\phi_{s}=\Xi_{s}^{-1}$
with $\Xi_{s}$ given in \eqref{eq:-9}.  By \eqref{eq:-28}, it
holds that 
\begin{align}
\mathbb{P}\{\ln f-\ln\|f\|_{p}\ge r\} & \le\rme^{nq\int_{p}^{q}\phi_{s}(\beta(s))s^{-2}\d s-rq},\label{eq:-28-1-1-1}
\end{align}
 where $\beta(s)$ is given in \eqref{eq:-10}. Optimizing the bound
over all $q\in[p,2]$ yields the desired bound. 

\section*{Acknowledgement}

The authors would like to thank anonymous reviewers for improving
the quality of this paper. Specifically, the authors thank one of
reviewers for pointing out the potential connections between the R\'enyi--Sobolev
inequalities and concentration inequalities, which motivates the investigation
in Section \ref{subsec:Connections-to-Concentration}. The authors 
also thank  another reviewer for pointing out the non-convexity of the
log-Sobolev function $\Xi_{q}$. 

\bibliographystyle{unsrt}
\bibliography{ref}

\begin{thebibliography}{10}

\bibitem{bakry2004functional}
D.~Bakry.
\newblock Functional inequalities for {Markov} semigroups.
\newblock In {\em Probability Measures on Groups}, pages 91--147. Tata
  Institute of Fundamental Research, Mumbai, 2004.

\bibitem{bakry2013analysis}
D.~Bakry, I.~Gentil, and M.~Ledoux.
\newblock {\em Analysis and Geometry of Markov Diffusion Operators}, volume
  348.
\newblock Springer Science \& Business Media, 2013.

\bibitem{rudnicki2002markov}
R.~Rudnicki, M.~Pich\'or, and M.~Tyran-Kami\'nska.
\newblock {\em Markov Semigroups and Their Applications}, volume 597 of {\em
  Dynamics of Dissipation. Lecture Notes in Physics}.
\newblock Springer, Berlin, Heidelberg, 2002.

\bibitem{Ledoux_book}
M.~Ledoux.
\newblock {\em Concentration of Measure and Logarithmic Sobolev Inequalities},
  volume 1709 of {\em S\'eminaire de Probabilit\'es XXXIII. Lecture Notes in
  Mathematics}.
\newblock Springer, Berlin, Heidelberg, 2006.

\bibitem{gross1975logarithmic}
L.~Gross.
\newblock Logarithmic {Sobolev} inequalities.
\newblock {\em American Journal of Mathematics}, 97(4):1061--1083, 1975.

\bibitem{mossel2013reverse}
E.~Mossel, K.~Oleszkiewicz, and A.~Sen.
\newblock On reverse hypercontractivity.
\newblock {\em Geometric and Functional Analysis}, 23(3):1062--1097, 2013.

\bibitem{polyanskiy2019improved}
Y.~Polyanskiy and A.~Samorodnitsky.
\newblock Improved log-{S}obolev inequalities, hypercontractivity and
  uncertainty principle on the hypercube.
\newblock {\em Journal of Functional Analysis}, 277(11):108280, 2019.

\bibitem{gu2023non}
Y.~Gu and Y.~Polyanskiy.
\newblock Non-linear {log-Sobolev} inequalities for the {Potts} semigroup and
  applications to reconstruction problems.
\newblock {\em Communications in Mathematical Physics}, 404(2):769--831, 2023.

\bibitem{yu2021strong_article}
L.~Yu.
\newblock Strong {Brascamp--Lieb} inequalities.
\newblock {\em ArXiv e-prints, arXiv:2102.06935}, 2021.

\bibitem{samorodnitsky2008modified_article}
A.~Samorodnitsky.
\newblock A modified logarithmic {Sobolev} inequality for the {Hamming} cube
  and some applications.
\newblock {\em ArXiv e-prints, arXiv:0807.1679}, 2008.

\bibitem{kirshner2021moment}
N.~Kirshner and A.~Samorodnitsky.
\newblock A moment ratio bound for polynomials and some extremal properties of
  {Krawchouk} polynomials and {Hamming} spheres.
\newblock {\em IEEE Transactions on Information Theory}, 67(6):3509--3541,
  2021.

\bibitem{yu2024graphs}
L.~Yu, V.~Anantharam, and J.~Chen.
\newblock Graphs of joint types, noninteractive simulation, and stronger
  hypercontractivity.
\newblock {\em IEEE Transactions on Information Theory}, 70(4):2287--2308,
  2024.

\bibitem{friedman2005generalized}
J.~Friedman and J.-P. Tillich.
\newblock Generalized alon--boppana theorems and error-correcting codes.
\newblock {\em SIAM Journal on Discrete Mathematics}, 19(3):700--718, 2005.

\bibitem{bollobas2018eigenvalues}
B.~Bollob{\'a}s, J.~Lee, and S.~Letzter.
\newblock Eigenvalues of subgraphs of the cube.
\newblock {\em European Journal of Combinatorics}, 70:125--148, 2018.

\bibitem{chafai2004entropies}
D.~Chafa{\"\i}.
\newblock Entropies, convexity, and functional inequalities, on
  {$\Phi$-entropies and $\Phi$-Sobolev} inequalities.
\newblock {\em Journal of Mathematics of Kyoto University}, 44(2):325--363,
  2004.

\bibitem{boucheron2013concentration}
S.~Boucheron, G.~Lugosi, and P.~Massart.
\newblock {\em Concentration inequalities: A nonasymptotic theory of
  independence}.
\newblock Oxford university press, 2013.

\bibitem{raginsky2016strong}
M.~Raginsky.
\newblock Strong data processing inequalities and {$ \Phi $-Sobolev}
  inequalities for discrete channels.
\newblock {\em IEEE Transactions on Information Theory}, 62(6):3355--3389,
  2016.

\bibitem{ledoux2001concentration}
M.~Ledoux.
\newblock {\em The concentration of measure phenomenon}.
\newblock Number~89. American Mathematical Soc., 2001.

\bibitem{RagSason}
M.~Raginsky and I.~Sason.
\newblock {\em Concentration of Measure Inequalities in Information Theory,
  Communications and Coding}, volume~10 of {\em Foundations and Trends in
  Communications and Information Theory}.
\newblock Now Publishers Inc, 2013.

\bibitem{sudakov1978extremal}
V.~N. Sudakov and B.~S. Tsirel'son.
\newblock Extremal properties of half-spaces for spherically invariant
  measures.
\newblock {\em Journal of Soviet Mathematics}, 9(1):9--18, 1978.

\bibitem{borell1975brunn}
C.~Borell.
\newblock The {Brunn-Minkowski} inequality in {Gauss} space.
\newblock {\em Inventiones mathematicae}, 30(2):207--216, 1975.

\bibitem{bollobas1986combinatorics}
B.~Bollob{\'a}s.
\newblock {\em Combinatorics: set systems, hypergraphs, families of vectors,
  and combinatorial probability}.
\newblock Cambridge University Press, 1986.

\bibitem{alon1998asymptotic}
N.~Alon, R.~Boppana, and J.~Spencer.
\newblock An asymptotic isoperimetric inequality.
\newblock {\em Geometric \& Functional Analysis}, 8(3):411--436, 1998.

\bibitem{gozlan2005principe}
N.~Gozlan.
\newblock {\em Principe conditionnel de Gibbs pour des contraintes fines
  approch{\'e}es et in{\'e}galit{\'e}s de transport}.
\newblock PhD thesis, Universit{\'e} de Nanterre-Paris X, 2005.

\bibitem{gozlan2007large}
N.~Gozlan and C.~L{\'e}onard.
\newblock A large deviation approach to some transportation cost inequalities.
\newblock {\em Probability Theory and Related Fields}, 139(1-2):235--283, 2007.

\bibitem{yu2024exact}
L.~Yu.
\newblock Exact exponents for concentration and isoperimetry in product
  {Polish} spaces.
\newblock {\em IEEE Transactions on Information Theory}, 70(8):5427--5452,
  2024.

\bibitem{liu2023unsolved}
L.~Liu and B.~Ning.
\newblock Unsolved problems in spectral graph theory.
\newblock {\em arXiv preprint arXiv:2305.10290}, 2023.

\bibitem{navon2009linear}
M.~Navon and A.~Samorodnitsky.
\newblock Linear programming bounds for codes via a covering argument.
\newblock {\em Discrete \& Computational Geometry}, 41:199--207, 2009.

\bibitem{elgamal}
A.~{El~Gamal} and Y.-H. Kim.
\newblock {\em Network Information Theory}.
\newblock Cambridge University Press, Cambridge, U.K., 2012.

\end{thebibliography}

\begin{IEEEbiographynophoto}{Lei Yu} (Member, IEEE)  received the B.E. and Ph.D. degrees in electronic  
engineering from the University of Science and Technology of China (USTC)  
in 2010 and 2015, respectively. From 2015 to 2020, he worked as a 
Post-Doctoral Researcher at the USTC, National University of Singapore, and  
University of California at Berkeley. He is currently an Associate  
Professor at the School of Statistics and Data Science, LPMC, KLMDASR,  
and LEBPS, Nankai University, China. Since 2024, he has served as  Associate Editor of  the IEEE Transactions on Information Theory. His research interests lie in the  
intersection of probability theory, information theory, and combinatorics.
\end{IEEEbiographynophoto}

\begin{IEEEbiographynophoto}{Hao Wu} received the B.S. degree in Mathematics from Sichuan University  in 2017 and the M.S. degree in Mathematics from Nankai University in 2021. He is currently a Ph.D. student at the School of Statistics and Data Science, Nankai University, China. His research interests include probability theory and information theory.  \end{IEEEbiographynophoto}
\end{document}